\documentclass[reqno, 11pt]{amsart}
\usepackage{amsmath}
\usepackage{amssymb}
\usepackage{graphicx}
\usepackage{amsfonts}
\usepackage[margin=0.93in, letterpaper ]{geometry}
\usepackage{lineno}
\usepackage[onehalfspacing]{setspace}

\newtheorem{theorem}{Theorem}
\theoremstyle{plain}

\newtheorem{corollary}[theorem]{Corollary}

\newtheorem{definition}[theorem]{Definition}

\newtheorem{lemma}[theorem]{Lemma}

\newtheorem{proposition}[theorem]{Proposition}

\newtheorem{remark}[theorem]{Remark}

\newtheorem{fact}[theorem]{Fact}
\numberwithin{equation}{section}
\numberwithin{theorem}{subsection}
\input{tcilatex}

\begin{document}
\title{\textsf{A look at Representations of }$SL_{2}(\mathbb{F}_{q})$\textsf{%
\ through the Lens of Size}}
\author{\textit{Shamgar Gurevich}}
\address{\textit{Department of Mathematics, University of Wisconsin,
Madison, WI 53706, USA.}}
\email{shamgar@math.wisc.edu}
\author{\textit{Roger Howe}}
\address{\textit{Department of Mathematics, Yale University, New Haven, CT
06520, USA.}}
\email{roger.howe@yale.edu}
\date{\textit{New Haven - CT, Saturday June 1, 2018.}}

\begin{abstract}
How to study a nice function on the real line? The physically motivated
Fourier theory technique of harmonic analysis is to expand the function in
the basis of exponentials and study the meaningful terms in the expansion.
Now, suppose the function lives on a finite non-commutative group $G$, and
is invariant under conjugation. There is a well-known analog of Fourier
analysis, using the irreducible characters of $G$. This can be applied to
many functions that express interesting properties of $G$. To study these
functions one wants to know how the different characters contribute to the
sum?

In this note we describe the $G=SL_{2}(\mathbb{F}_{q})$ case of the theory
we have been developing in recent years which attempts to give a fairly
general answer to the above question for finite classical groups.

The irreducible representations of $SL_{2}(\mathbb{F}_{q})$ are
\textquotedblleft well known\textquotedblright\ for a very long time \cite%
{Frobenius1896, Jordan1907, Schur1907} and are a prototype example in many
introductory courses on the subject. We are happy that we can say something
new about them. In particular, it turns out that the representations that
were considered as \textquotedblleft anomalous"\ in the "old" point of view
(known as the "philosophy of cusp forms") are the building blocks of the
current approach.
\end{abstract}

\maketitle
\thanks{ \ \ \ \ \ \ \ \ \ \ \ \ \ \ \ \ \ \ \ \ \ \ \ \ \ \ \ \ \ \ \ \ \ \
\ \ \textit{\ \ To Joe Wolf mazal tov ad me'ah v'esrim}}

\section{\textbf{Introduction}}

An important invariant attached to any finite group $G$ is its \textit{%
representation ring} \cite{Zelevinsky81} 
\begin{equation*}
R(G)=%
%TCIMACRO{\U{2124} }%
%BeginExpansion
\mathbb{Z}
%EndExpansion
\lbrack \widehat{G}],
\end{equation*}%
i.e., the ring generated from the set $\widehat{G}$ of isomorphism classes
of irreducible representations (irreps) using the operations of addition and
multiplication given, respectively, by direct sum $\oplus $ and tensor
product $\otimes .$

The main goal of this note is to advertise, using the special linear group
over the finite field $\mathbb{F}_{q}$ with $q$ odd elements 
\begin{equation*}
SL_{2}(\mathbb{F}_{q})=\left\{ 
\begin{pmatrix}
a & b \\ 
c & d%
\end{pmatrix}%
;\text{ }ad-bc=1,\text{ \ }a,b,c,d\in \mathbb{F}_{q}\right\} ,
\end{equation*}%
our recent discovery \cite{Gurevich-Howe15, Gurevich-Howe17, Gurevich17,
Howe17-1, Howe17-2} that in the case that $G$ is a finite classical group
the ring $R(G)$ has a natural "tensor rank" filtration that encodes the
analytic properties of the representations. In particular, moving to the
associated graded pieces, each member $\rho $ of the unitary dual $\widehat{G%
}$ gets a well defined non-negative integer called \textbf{tensor rank} and
denoted $rank_{\otimes }(\rho ).$ This integer seems---see Figure \ref%
{cr-u-sl2-101} for illustration\footnote{%
The numerics in this note were generated using the computer algebra system
Magma.} in the case of $G=SL_{2}(\mathbb{F}_{q})$---to be intimately related
to analytic properties of $\rho $ such as its dimension and size of
character values.%
%TCIMACRO{%
%\FRAME{ftbpFU}{5.0834in}{3.8268in}{0pt}{\Qcb{$\log _{1/\sqrt{101}}$-scale of character ratios at $u=\tbinom{1\text{ \ }1}{0\text{ \ }1}$ vs. tensor rank for irreps of $SL_{2}(\mathbb{F}_{101})$.}}{\Qlb{cr-u-sl2-101}}{cr-u-sl2-101.png}{%
%\special{language "Scientific Word";type "GRAPHIC";maintain-aspect-ratio TRUE;display "USEDEF";valid_file "F";width 5.0834in;height 3.8268in;depth 0pt;original-width 5.028in;original-height 3.7775in;cropleft "0";croptop "1";cropright "1";cropbottom "0";filename '//vmware-host/Shared Folders/shamgar On My Mac/Dropbox/Collaborators/Howe/Project (1)/Papers/SL2/Arxiv/graphics/CR-u-SL2-101.png';file-properties "XNPEU";}}}%
%BeginExpansion
\begin{figure}[ptb]\centering
\includegraphics
%[natheight=3.7775in, natwidth=5.028in, height=3.8268in, width=5.0834in]
{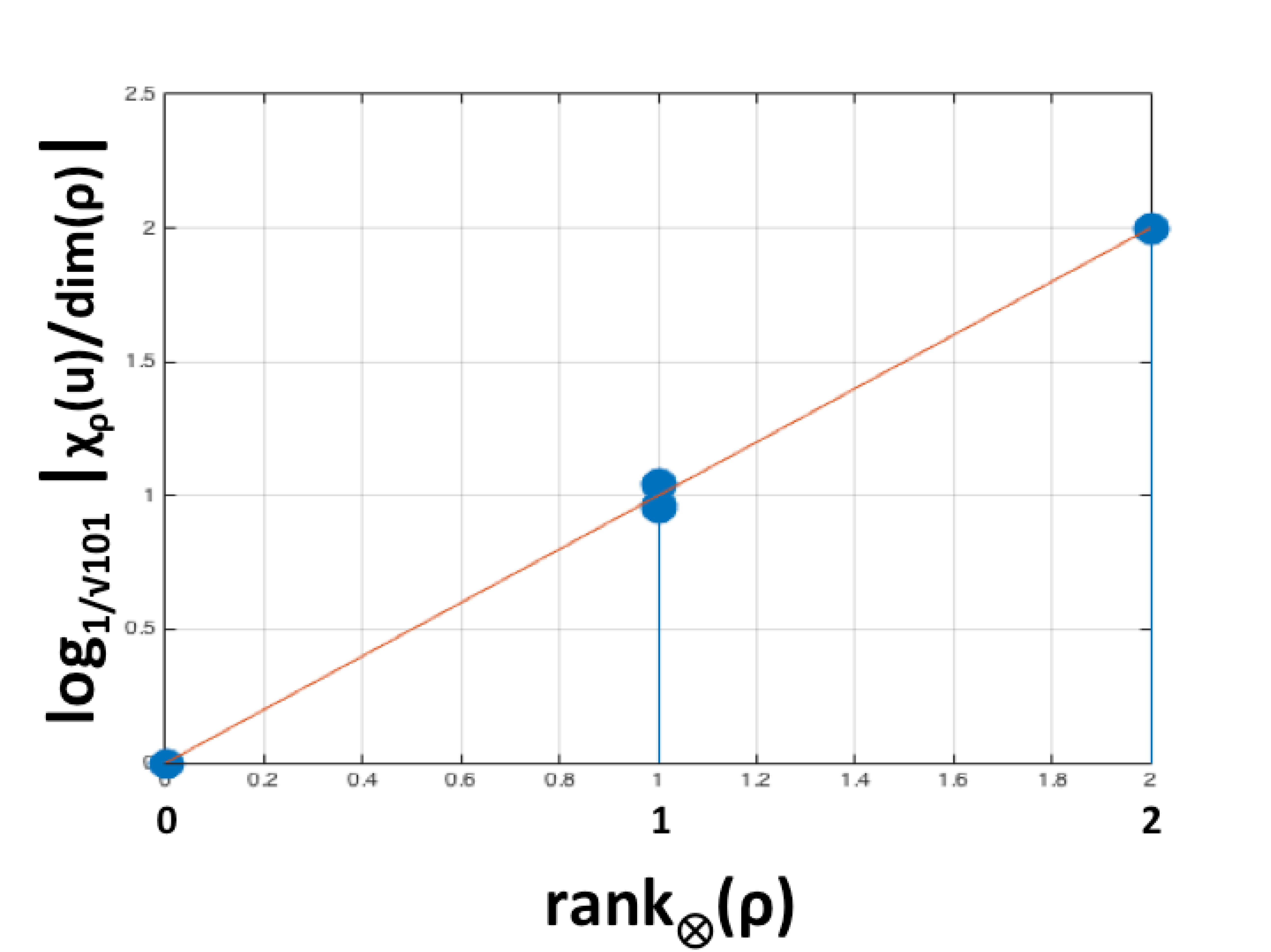}%
\caption{$\log _{1/\protect\sqrt{101}}$-scale of character ratios at $u=%
\tbinom{1\text{ \ }1}{0\text{ \ }1}$ vs. tensor rank for irreps of $SL_{2}(%
\mathbb{F}_{101})$.}\label{cr-u-sl2-101}%
\end{figure}%
%EndExpansion

In the rest of the introduction we discuss two aspects of the tensor rank
filtration. The first is the fact that this is the first point of view that
makes the degenerate discrete and principal series representations of $%
SL_{2}(\mathbb{F}_{q})$ seem significant rather than anomalous. The second
aspect we discuss, through a specific example for $SL_{2}(\mathbb{F}_{q}),$
is the potential usefulness of the tensor rank filtration to the harmonic
analysis of the finite classical groups.

\subsection{\textbf{Degenerate Discrete and Principal Series Representations 
\label{DD-DP}}}

Let $G=\mathbf{G}(\mathbb{F}_{q})$ be the group of $\mathbb{F}_{q}$-rational
points of a reductive algebraic group defined over $\mathbb{F}_{q}.$ The
Gelfand--Harish-Chandra "philosophy of cusp forms" \cite{Gelfand62,
Harish-Chandra70} asserts that each member $\rho $ of $\widehat{G}$ can be
realized, in some effective manner, inside the induction $Ind_{P}^{G}(\sigma
)$ from a cuspidal representation (i.e., one that is not induced from a
smaller parabolic subgroup) of a parabolic subgroup $P$ of $G$ 
\begin{equation*}
\widehat{G}=\left\{ 
\begin{array}{c}
\rho <Ind_{P}^{G}(\sigma )\text{, } \\ 
G\supset P\text{ - parabolic,} \\ 
\sigma \text{ cuspidal repn of }P%
\end{array}%
\right\} .
\end{equation*}%
The above philosophy is a central one and leads to important developments in
representation theory of reductive groups over local and finite fields, in
particular the work of Deligne--Lusztig \cite{Deligne-Lusztig76} on the
construction of representations of finite reductive groups and Lusztig's
striking achievement: The classification \cite{Lusztig84} of these
representations. 
%TCIMACRO{%
%\FRAME{fhFU}{5.8963in}{1.8127in}{0pt}{\Qcb{Rough description of $\widehat{GL_{2}(\mathbb{F}_{q})}$ and $\widehat{SL_{2}(\mathbb{F}_{q})}.$}}{\Qlb{philosophy}}{philosophy.bmp}{%
%\special{language "Scientific Word";type "GRAPHIC";maintain-aspect-ratio TRUE;display "USEDEF";valid_file "F";width 5.8963in;height 1.8127in;depth 0pt;original-width 11.0558in;original-height 3.5414in;cropleft "0";croptop "1";cropright "1";cropbottom "0";filename '../Gurevich-Howe's paper/Philosophy.bmp';file-properties "XNPEU";}}}%
%BeginExpansion
\begin{figure}[h]\centering
\includegraphics
%[natheight=3.5414in, natwidth=11.0558in, height=1.8127in, width=5.8963in]
{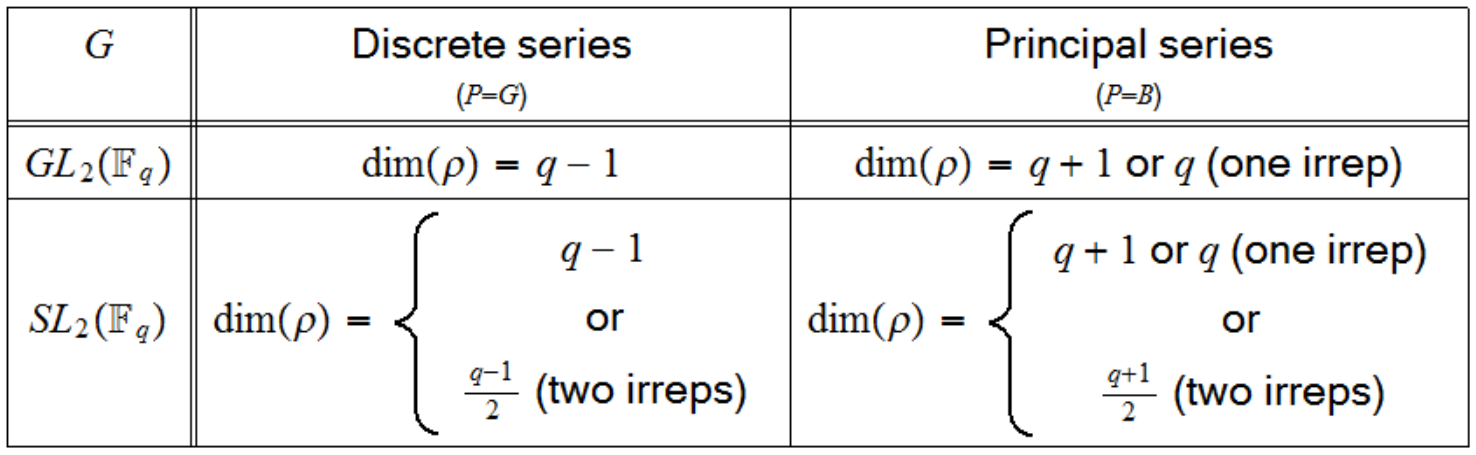}%
\caption{Rough description of $\protect\widehat{GL_{2}(\mathbb{F}_{q})}$ and 
$\protect\widehat{SL_{2}(\mathbb{F}_{q})}.$}\label{philosophy}%
\end{figure}%
%EndExpansion

Let us see how this philosophy manifests itself in the cases when $G$ is $%
GL_{2}(\mathbb{F}_{q})$ or $SL_{2}(\mathbb{F}_{q})$-see Figure \ref%
{philosophy} (let us ignore for this discussion the one dimensional
representations). In both cases around $50\%$ of the irreps are cuspidals,
i.e., $P=G,$ and form the so called discrete series, and $50\%$ are the so
called principal series, i.e., induced from the Borel subgroup $P=B$ of
upper triangular matrices in $G.$ Moreover, it is well known, and not so
difficult to show \cite{Fulton-Harris91}, that the discrete and principal
series representations of $SL_{2}(\mathbb{F}_{q})$ can be obtained by
restriction from the corresponding irreps of $GL_{2}(\mathbb{F}_{q}).$ These
restrictions are typically irreducible, with only two exceptions: One
discrete series representation of dimension $q-1$ and one principal series
representation of dimension $q+1$ split into two pieces. Hence one has---see
Figure \ref{philosophy}---two discrete and two principal series
representations of dimensions $\frac{q-1}{2}$ and $\frac{q+1}{2}$,
respectively. These representations are called in the literature "degenerate
discrete series" and "degenerate principal series" and are usually treated
as anomalous.

In this note we show that the tensor rank filtration proposes a quantitative
measurement that, in particular, demonstrates that the degenerate
representations are significant rather then anomalous. In fact, in some
formal sense these irreps are the "atoms" of the all theory of "size" of
representations we are describing.

\subsection{\textbf{Harmonic Analysis on a Finite Group}}

Suppose we want to study a class function $\mathcal{N}$, defined on a finite
group $G$: 
\begin{equation*}
\mathcal{N}:G\rightarrow 
%TCIMACRO{\U{2102} }%
%BeginExpansion
\mathbb{C}
%EndExpansion
,\text{ \ }\mathcal{N(}hgh^{-1})=\mathcal{N}(g),
\end{equation*}%
for every $g,h\in G.$

Many interesting examples of $\mathcal{N}$ present similar manner discussed
below.

\subsubsection{\textbf{The problem of harmonic analysis on a finite group}}

Class functions can, in principle, be investigated using the representation
theory of $G$: it is a basic fact \cite{Serre77} that one can expand them as
a linear combination of the irreducible characters of $G$. This is the
harmonic analysis approach for studying class functions on $G.$ In
particular, this technique can be applied to many functions that express
interesting properties of $G$ \cite{Diaconis88, Guralnick-Malle14,
Liebeck17, Liebeck-O'Brien-Shalev-Tiep10, Malle14, Shalev07, Shalev17}. In
that latter case the expansion is in many cases of the form 
\begin{equation*}
\mathcal{N=}\dsum\limits_{\rho \in \widehat{G}}"\frac{\chi _{\rho }}{\dim
(\rho )}",
\end{equation*}%
where $"\frac{\chi _{\rho }}{\dim (\rho )}"$ stands for some relevant
expression in terms of the \underline{character ratios} $\frac{\chi _{\rho }%
}{\dim (\rho )}$ of the irreps of $G.$

The above discussion suggests the following:\medskip

\textbf{Problem (Main problem of harmonic analysis on }$G$\textbf{). }%
Estimate the character ratios%
\begin{equation*}
\frac{\chi _{\rho }(g)}{\dim (\rho )},\text{ \ }\rho \in \widehat{G},\text{ }%
g\in G,
\end{equation*}%
and possible relations among them.

\subsubsection{\textbf{Example: The commutator mapping on }$SL_{2}(\mathbb{F}%
_{q})$\label{E-CM}}

Consider the group $G=SL_{2}(\mathbb{F}_{q})$ and the commutator mapping 
\begin{equation}
\left[ ,\right] :G\times G\rightarrow G,\ \left[ x,y\right] =xyx^{-1}y^{-1}.
\label{CM}
\end{equation}%
For an element $g\in G$ let us denote by $[,]_{g}$ the set $\{(x,y)\in
G\times G;$ $[x,y]=g\}$ called the \textit{fiber over }$g$ of the commutator
mapping. We are interested in the distribution of the number of elements $%
\#([,]_{g})$ when $g$ runs in the set $G\smallsetminus \{\pm I\}.$ It is
natural \cite{Shalev07} to make a suitable normalization and study the class
function 
\begin{equation}
\mathcal{N}(g)=\#(\left[ ,\right] _{g})/\#G.  \label{N}
\end{equation}%
To see how this function behaves over the various conjugacy classes of $G$,
let us recall \cite{Fulton-Harris91} that $G$ has $q+4$ conjugacy classes
with standard representatives 
\begin{equation}
I,\text{ \ }-I,\text{ \ }%
\begin{pmatrix}
\pm 1 & 1,\varepsilon \\ 
0 & \pm 1%
\end{pmatrix}%
,\text{ \ }%
\begin{pmatrix}
a & 0 \\ 
0 & a^{-1}%
\end{pmatrix}%
,\text{ \ }%
\begin{pmatrix}
x & \varepsilon y \\ 
y & x%
\end{pmatrix}%
,  \label{CC}
\end{equation}%
where $\varepsilon $ is a non-square\footnote{%
For the rest of this note $\varepsilon $ stands for a non-square element of $%
\mathbb{F}_{q}.$}, $a\neq \pm 1,$ $y\neq 0,$ and $x^{2}-\varepsilon y^{2}=1.$
%TCIMACRO{%
%\FRAME{fhFU}{5.7839in}{1.5281in}{0pt}{\Qcb{Values of $\mathcal{N}(g)$ for various $g,q.$}}{\Qlb{n-by-g}}{n-by-g.bmp}{%
%\special{language "Scientific Word";type "GRAPHIC";maintain-aspect-ratio TRUE;display "USEDEF";valid_file "F";width 5.7839in;height 1.5281in;depth 0pt;original-width 9.2915in;original-height 3.4861in;cropleft "0";croptop "1";cropright "1";cropbottom "0";filename '../../../Lectures/Size of Repn of SL(2,q)/graphics/N-by-G.bmp';file-properties "XNPEU";}}}%
%BeginExpansion
\begin{figure}[h]\centering
\includegraphics
%[natheight=3.4861in, natwidth=9.2915in, height=1.5281in, width=5.7839in]
{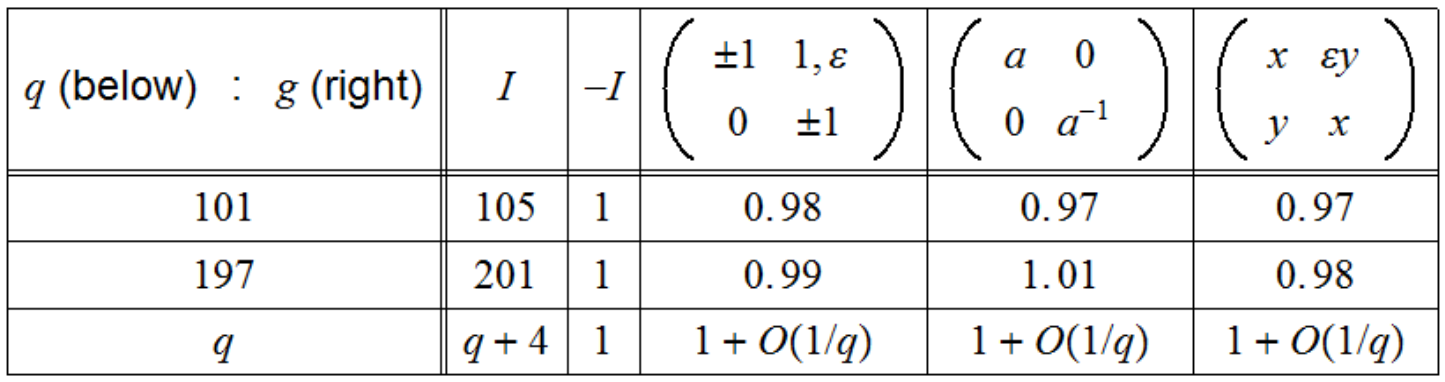}%
\caption{Values of $\mathcal{N}(g)$ for various $g,q.$}\label{n-by-g}%
\end{figure}%
%EndExpansion

Using numerics---see Figure \ref{n-by-g} for illustration---one can suspect
the following to be true:\medskip

\textbf{Theorem (Uniformity of the commutator map). }We have\footnote{%
For numerical sequences $a(q),b(q),$ we write $a(q)=O(b(q))$ if there is $C$
s.t. $\left\vert a(q)\right\vert \leq C\left\vert b(q)\right\vert $ for all
sufficiently large $q.$} 
\begin{equation}
\mathcal{N}(g)=1+O(1/q),\ g\neq \pm I.  \label{UThm}
\end{equation}%
$\medskip $

The relevant formula enabling the representation theoretic study of the
function $\mathcal{N}$ is due to Frobenius \cite{Frobenius1896} 
\begin{equation}
\mathcal{N}(g)=\sum\limits_{\rho \in \widehat{G}}\frac{\chi _{\rho }(g)}{%
\dim (\rho )}.  \label{FSum}
\end{equation}%
Looking back on (\ref{UThm}), and taking into account the trivial
representation in (\ref{FSum}), we see that to verify the uniformity
statement above we need to\medskip

\begin{description}
\item[\textbf{Goal}] Explore the cancellation%
\begin{equation}
\sum\limits_{1\neq \rho \in \widehat{G}}\frac{\chi _{\rho }(g)}{\dim (\rho )}%
=O(1/q),\text{ \ }g\neq \pm I.  \label{Goal-C}
\end{equation}
\end{description}

In fact, it is possible (see \cite{Garion-Shalev09}) to verify (\ref{Goal-C}%
) invoking available character tables of $G$ \cite{Fulton-Harris91}. There
are two drawbacks to this approach. First, it will prevent us from telling
the story we are trying to make. Second, and more seriously, it is not
really feasible to have explicit character tables for more complicated
groups, e.g., the high rank symplectic groups. This is a family that
contains $SL_{2}(\mathbb{F}_{q})=Sp_{2}(\mathbb{F}_{q})$ as a member and for
which one might expect (see the numerics in \cite{Gurevich-Howe15}) similar
behavior of the relevant Frobenius sum (\ref{Goal-C}). So we would like to
have a different approach to achieve the goal (\ref{Goal-C}). In this note
we demonstrate, in the case of $SL_{2}(\mathbb{F}_{q}),$ how the tensor rank
filtration might assists with that task.\medskip

\textbf{Acknowledgements. }This note was written in the winter of 2017-18
while S.G. was visiting the Mathematics Departments at Weizmann Institute
and Yale University, and the College of Education at Texas AM University,
and he would like to thank these institutions, personally H. Naor and O.
Sarig at Weizmann, D. Altschuler, C. Villano, and I. Frenkel at Yale, and R.
Howe at TAMU. We also want to thank S. Goldstein and J. Cannon for their
help with numerical aspects of the project, part of it is reported here. \ 

\section{\textbf{Character Ratios and Tensor Rank\label{CR-TR}}}

A priori, it is possible that the sum (\ref{Goal-C}) is small as it is
because each of the (around) $q$ terms there is so small. This is not the
case.

\subsection{\textbf{Character Ratios}}

\begin{figure}[h]\centering
\includegraphics
%[natheight=2.6671in, natwidth=10.8612in, height=1.8498in, width=4.0785in]
{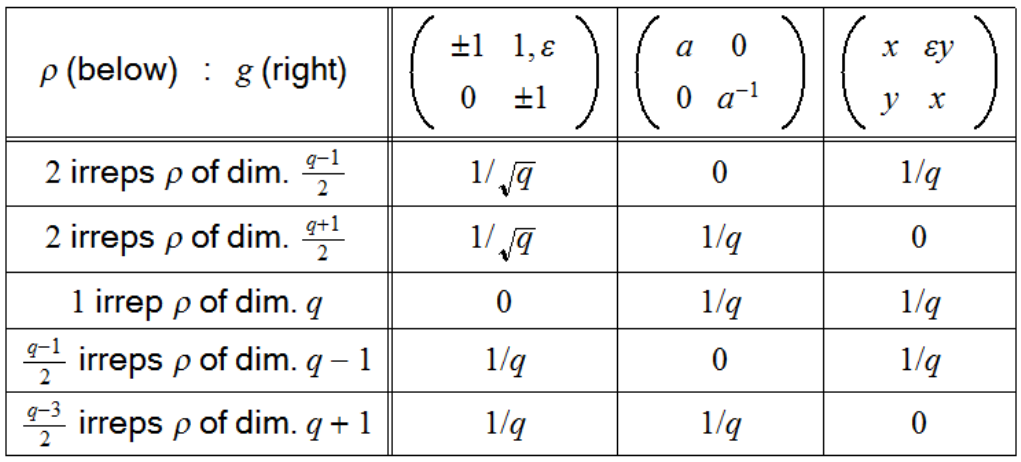}%
\caption{Order of magnitude of $\protect\chi _{\protect\rho }(g)/\dim (%
\protect\rho )$ for $\protect\rho $'s irreps of $G=SL_{2}(\mathbb{F}_{q})$. }%
\label{cr-sl2-q}%
\end{figure}

A quick numerical experiment with Magma tells us (see Figure \ref{cr-sl2-q}
for the outcome) that indeed the sum (\ref{Goal-C}) is small due to
cancellations. For example, 
%TCIMACRO{%
%\FRAME{fhFU}{4.0785in}{1.8498in}{0pt}{\Qcb{Order of magnitude of $\chi _{\rho }(g)/\dim (\rho )$ for $\rho $'s irreps of $G=SL_{2}(\mathbb{F}_{q})$. }}{\Qlb{cr-sl2-q}}{cr-sl2-q.bmp}{%
%\special{language "Scientific Word";type "GRAPHIC";maintain-aspect-ratio TRUE;display "USEDEF";valid_file "F";width 4.0785in;height 1.8498in;depth 0pt;original-width 10.8612in;original-height 2.6671in;cropleft "0";croptop "1";cropright "1";cropbottom "0";filename '../Gurevich-Howe's paper/CR-SL2-q.bmp';file-properties "XNPEU";}}}%
%BeginExpansion
%EndExpansion
in the last two rows of the second column of Figure \ref{cr-sl2-q} there are
around $q$ character ratios each one of order of magnitude $1/q.$ Moreover,
looking again on the second and third rows of the second column, one hope
for two "sub-cancellations" in the sum (\ref{Goal-C}), one between the
character ratios of the irreps of dimensions $\frac{q\pm 1}{2},$ and another
one between the character ratios of the other irreps.

We will show that indeed this what is going on.

\subsection{\textbf{Tensor Rank of a Representation}}

The different orders of magnitude of the character ratios discussed
above---see also Figure \ref{cr-u-sl2-101} for illustration---is a
manifestation of a particular invariant---called tensor rank---that can be
attached to any representation of $G$ and should be regarded in a formal
sense as its "size".

Let $k$ be a non-negative integer and consider the collection $%
F^{k}=F^{k}R(G)$ of representations in $R(G)$ that appear inside sums (with
integer coefficients) of $\ell $-fold, $\ell \leq k$, tensor product of the
four irreps of $G$ of dimensions $\frac{q\pm 1}{2}.$ One has $F^{k}\subset
F^{k+1}$, $F^{i}F^{j}\subset F^{i+j}$ and $\cup _{k}F^{k}=R(G)$, so the $%
F^{k}$'s form a filtration of $R(G)$ that we call \textit{tensor rank
filtration. }In fact,

\begin{proposition}
\label{F2R}\textit{\ } We have $F^{2}=R(G).$
\end{proposition}

Proposition \ref{F2R} is verified using explicit computations in Section \ref%
{Exp-etaC}---see Subsection \ref{Summary}.\smallskip

Now we can pass to the \textit{associated graded }pieces $%
Gr^{k}=F^{k}/F^{k-1}$ and in particular obtain a well defined notion of
tensor rank for representations \cite{Gurevich-Howe17}. Concretely,

\begin{definition}[\textbf{Tensor rank}]
Let $\rho $ be a representation of $G.$ We say that $\rho $ has \underline{%
tensor rank} $k,$ denoted $rank_{\otimes }(\rho )=k,$ if \ $\rho \in
F^{k}\smallsetminus F^{k-1}.$
\end{definition}

Let us denote by $\widehat{G}_{k}$ the family of irreps of tensor rank $k.$
Combining Proposition \ref{F2R} and the information appearing in Figure \ref%
{philosophy} we see that $\widehat{G}_{0}$ consists of the trivial
representation, $\widehat{G}_{1}$ is the collection of irreps of dimensions $%
\frac{q\pm 1}{2}$, and finally $\widehat{G}_{2}$ is the family of irreps of
dimensions $q\pm 1$ and $q.$

Going back to the problem of exploring the cancellation in the sum (\ref%
{Goal-C}) our idea is to split it over the tensor ranks 
\begin{equation*}
\sum\limits_{1\neq \rho \in \widehat{G}}\frac{\chi _{\rho }(g)}{\dim (\rho )}%
=\sum\limits_{\rho \in \widehat{G}_{1}}\frac{\chi _{\rho }(g)}{\dim (\rho )}%
+\sum\limits_{\rho \in \widehat{G}_{2}}\frac{\chi _{\rho }(g)}{\dim (\rho )},
\end{equation*}%
and to witness the cancellation in each partial sum. Of course, for this we
will need more information on the irreps of each given tensor rank.

\section{\textbf{The Heisenberg and Oscillator Representations\label{HOR}}}

Where do the tensor rank one irreps of $G=SL_{2}(\mathbb{F}_{q})$ come from?
A. Weil shed some light on this in \cite{Weil64}. They can be found by
considering the Heisenberg\textit{\ }group.

\subsection{\textbf{The Heisenberg Group}}

We consider the vector space $W=\mathbb{F}_{q}^{2n}$ of column vectors of
length $2n$ with entries in $\mathbb{F}_{q}$. We equip $W$ with the
non-degenerate skew-symmetric bilinear form 
\begin{equation*}
\left\langle w,w^{\prime }\right\rangle =w^{t}%
\begin{pmatrix}
0 & I \\ 
-I & 0%
\end{pmatrix}%
w^{\prime },
\end{equation*}%
where $w^{t}$ is the row vector transpose of $w$, and $I$ is the $n\times n$
identity matrix. The pair $(W,\left\langle ,\right\rangle )$ is called
symplectic vector space.

The \textit{Heisenberg }group attached to $(W,\left\langle ,\right\rangle )$
is a two-step nilpotent group that can be realized by the set%
\begin{equation*}
H=W\times \mathbb{F}_{q},
\end{equation*}%
with the group law%
\begin{equation*}
\left( w,z\right) \cdot (w^{\prime },z^{\prime })=(w+w^{\prime },z+z^{\prime
}+\frac{1}{2}\left\langle w,w^{\prime }\right\rangle ).
\end{equation*}%
In particular, the center $Z$ of the Heisenberg group 
\begin{equation*}
Z=\{(0,z);\text{ }z\in \mathbb{F}_{q}\}=\mathbb{F}_{q},
\end{equation*}%
is equal to its commutator subgroup.

\subsection{\textbf{Representations of the Heisenberg Group\label{RHG}}}

We would like to describe the irreps of the Heisenberg group.

Take an irreducible representation $\pi $ of $H$. Then, by Schur's lemma 
\cite{Serre77}, the center $Z$ will act by scalars%
\begin{equation*}
\pi (z)=\psi _{\pi }(z)I,\text{ \ }z\in Z,
\end{equation*}%
where $I$ is the the identity operator on the representation space of $\pi $%
, and $\psi _{\pi }\in \widehat{Z}$ is a character of $Z,$ called the 
\textit{central character }of $\pi $. If $\psi _{\pi }=1,$ then $\pi $
factors through $H/Z\simeq W$, which is abelian, so $\pi $ is itself a
character of $W.$ The case of non-trivial central character is described by
the following celebrated theorem \cite{Mackey49}:

\begin{theorem}[\textbf{Stone--von Neumann--Mackey}]
Up to equivalence, there is a unique irreducible representation $\pi _{\psi
} $ with given non-trivial central character $\psi $ in $\widehat{Z}%
\smallsetminus \{1\}.$
\end{theorem}

We will call the (isomorphism class of the) representation $\pi _{\psi }$
the \textit{Heisenberg representation }associated to the central character $%
\psi .$

\begin{remark}[\textbf{Realization}]
There are many ways to realize (i.e., to write explicit formulas for) $\pi
_{\psi }$ \cite{Gerardin77, Gurevich-Hadani07, Gurevich-Hadani09, Howe73-1,
Weil64}. In particular, it can be constructed as induced representation from
any character extending $\psi $ to any maximal abelian subgroup of $H$. To
have a concrete one, note that the inverse image in $H$ of any maximal
subspace of $W$ on which $\left\langle ,\right\rangle $ is identically zero
(a.k.a. Lagrangian) will be a maximal abelian subgroup for which we can
naturally extend the character $\psi $. For example, consider the Lagrangian 
\begin{equation}
X=\left\{ 
\begin{pmatrix}
x_{1} \\ 
\vdots \\ 
x_{n} \\ 
0 \\ 
\vdots \\ 
0%
\end{pmatrix}%
\right\} \subset W,  \label{X}
\end{equation}%
and the associated maximal abelian subgroup $\widetilde{X}$ with character $%
\widetilde{\psi }$ on it, given by 
\begin{equation*}
\widetilde{X}=X\times \mathbb{F}_{q},\text{ \ }\widetilde{\psi }(x,z)=\psi
(z).
\end{equation*}%
Then we have the explicit realization of $\pi _{\psi },$ given by the action
of $H,$ by right translations, on the space 
\begin{equation}
Ind_{\widetilde{X}}^{H}(\widetilde{\psi })=\{f:H\rightarrow 
%TCIMACRO{\U{2102} }%
%BeginExpansion
\mathbb{C}
%EndExpansion
;\text{ \ }f(\widetilde{x}h)=\widetilde{\psi }(\widetilde{x})f(h)\text{, \ \ 
}\widetilde{x}\in \widetilde{X}\text{, }h\in H\}.  \label{Ind}
\end{equation}%
In particular, we have $\dim (\pi _{\psi })=q^{n}.$
\end{remark}

\subsection{\textbf{The Oscillator Representation}}

Consider the group $Sp(W)$ of automorphism of $W$ that preserve the form $%
\left\langle ,\right\rangle $. The action of $Sp(W)$ on $W$ lifts to an
action on $H$ by automorphisms leaving the center point-wise fixed. The
precise formula is $g(w,z)=(gw,z),$ \ $g\in Sp(W).$ It follows from the 
\textit{Stone--von Neumann--Mackey} theorem, that the induced action of $%
Sp(W)$ on the set $Irr(H)$ will leave fixed each isomorphism class $\pi
_{\psi },$ $\psi \in \widehat{Z}\smallsetminus \{1\}.$ This means that, if
we fix a vector space $\mathcal{H}_{\psi }$ realizing $\pi _{\psi },$ then
for each $g$ in $Sp(W)$ there is an operator $\omega _{\psi }(g)$ that acts
on space $\mathcal{H}_{\psi }$ and satisfies the equation 
\begin{equation}
\omega _{\psi }(g)\pi _{\psi }(h)\omega _{\psi }(g)^{-1}=\pi _{\psi }(g(h)).
\label{Egorov}
\end{equation}%
Note that, by Schur's lemma, the operator $\omega _{\psi }(g)$ is defined by
(\ref{Egorov}) up to scalar multiples. This implies that for any $%
g,g^{\prime }\in Sp(W)$ we have $\omega _{\psi }(g)\omega _{\psi }(g^{\prime
})=c(g,g^{\prime })\omega _{\psi }(gg^{\prime }),$ where $c(g,g^{\prime })$
is an appropriate complex number of absolute value $1$. It is well known
(see \cite{Gerardin77, Gurevich-Hadani07, Gurevich-Hadani09} for explicit
formulas) that over finite fields of odd characteristic the mapping $\omega
_{\psi }$ can be defined so that $c(g,g^{\prime })=1$ for every $g,g\in
Sp(W) $, i.e., $\omega _{\psi }$ defines a representation of $Sp(W).$ We
summarize:

\begin{theorem}[\textbf{Oscillator representation}]
\label{OR}There exists\footnote{%
This representation is unique except the case $n=2$ and $q=3,$ where there
is a canonical one \cite{Gurevich-Hadani07, Gurevich-Hadani09}.} a
representation 
\begin{equation}
\omega _{\psi }:Sp(W)\longrightarrow GL(\mathcal{H}),  \label{ome-psi}
\end{equation}%
that satisfies the identity (\ref{Egorov}).
\end{theorem}

We will call $\omega _{\psi }$ the \textit{oscillator representation}. This
is a name that was given to this representation in \cite{Howe73-1} due to
its origin in physics \cite{Segal60, Shale62}. Another popular name for $%
\omega _{\psi }$ is the \textit{Segal--Shale--Weil }or just\textit{\ Weil }%
representation, following the influential paper \cite{Weil64}.

\begin{remark}[\textbf{Schr\"{o}dinger model}]
\label{Schrodinger Model}We would like to have some useful formulas for $%
\omega _{\psi }$. We consider the Lagrangian decomposition $W=X\oplus Y$,
where $X$ is the first $n$ coordinates Lagrangian (\ref{X}) and $Y$ is the
second $n$ coordinates Lagrangian defined similarly. In terms of this
decomposition, the space (\ref{Ind}) is naturally identified with%
\begin{equation}
L^{2}(Y)\text{ \ }=\text{ \ functions on }Y.  \label{Space}
\end{equation}%
On the space (\ref{Space}) we realize the representation $\omega _{\psi }$.
This realization is sometime called the \textit{Schr\"{o}dinger model \cite%
{Howe80}. }In particular, in that model for every $f\in L^{2}(Y)$ we have 
\cite{Gerardin77, Gurevich-Howe15, Weil64}\bigskip

(A) $\left[ \omega _{\psi }%
\begin{pmatrix}
I & A \\ 
0 & I%
\end{pmatrix}%
f\right] (y)=\psi (\frac{1}{2}A(y,y))f(y),\medskip $

where $A:Y\rightarrow X$ is symmetric;\bigskip

(B) $\left[ \omega _{\psi }%
\begin{pmatrix}
0 & B \\ 
-B^{-1} & 0%
\end{pmatrix}%
f\right] (y)=\frac{1}{\gamma (B,\psi )}\tsum\limits_{y^{\prime }\in Y}\psi
(B(y,y^{\prime }))f(y^{\prime }),\medskip $

where $B:Y\widetilde{\rightarrow }X$ is symmetric, and $\gamma (B,\psi
)=\tsum\limits_{y\in Y}\psi (-\frac{1}{2}B(y,y))$ the quadratic Gauss sum 
\cite{Ireland-Rosen90};\bigskip

(C) $\left[ \omega _{\psi }%
\begin{pmatrix}
^{t}C^{-1} & 0 \\ 
0 & C%
\end{pmatrix}%
f\right] (y)=\QOVERD( ) {\det (C)}{q}f(C^{-1}y),\medskip $

where $C\in GL(Y),$ $^{t}C^{-1}\in GL(X)$ its transpose inverse, and $%
\QOVERD( ) {\cdot }{q}$ is the Legendre symbol\footnote{%
For $a\in \mathbb{F}_{q}^{\ast }$ the Legendre symbol $\QOVERD( ) {a}{q}=+1$
or $-1,$ according to $a$ being a square or not, respectively$.$}.
\end{remark}

It turns out that the isomorphism class of $\omega _{\psi }$ does depends on
the central character $\psi $ in $\widehat{Z}\smallsetminus \{1\}.$ However,
this dependence is weak. The following result \cite{Gerardin77, Howe73-1,
Howe73-2} indicates that there are only two possible oscillator
representations. For a character $\psi $ in $\widehat{Z}\smallsetminus \{1\}$
denote by $\psi _{a},$ $a\in \mathbb{F}_{q}^{\ast },$ the character $\psi
_{a}(z)=\psi (az).$

\begin{proposition}
\label{Oscillators}We have $\omega _{\psi }\simeq \omega _{\psi ^{\prime }}$
iff $\psi ^{\prime }=\psi _{a^{2}\text{ }}$ for some $a\in \mathbb{F}%
_{q}^{\ast }.$
\end{proposition}

\subsection{\textbf{The Representations of Tensor Rank One\label{RO}}}

The oscillator representations are slightly reducible. Indeed, for any two
representations $\varrho ,\varrho ^{\prime }$ of any group one has the 
\textit{intertwining number }$(\varrho ,\varrho ^{\prime })=\dim Hom(\varrho
,\varrho ^{\prime })$. The following is well known \cite{Gerardin77,
Howe73-1, Howe82}.

\begin{proposition}
\label{IO}We have $(\omega _{\psi },\omega _{\psi })=2.$
\end{proposition}

For the sake of completeness, we give a proof of Proposition \ref{IO} in
Appendix \ref{P-IO}.\medskip

The meaning of Proposition \ref{IO} is that $\omega _{\psi }$ has two
irreducible pieces. They can be computed explicitly as follows. The center $%
Z(Sp(W))=\{\pm I\}$ acts on the representation $\omega _{\psi }$, inducing
the decomposition into irreps 
\begin{equation}
\omega _{\psi }=\omega _{\psi ,1}\oplus \omega _{\psi ,sgn},  \label{O-Comp}
\end{equation}%
where $\omega _{\psi ,1}$ is the subspace of vectors on which $Z(Sp(W))$
acts trivially, and $\omega _{\psi ,sgn}$ is the subspace of vectors on
which $Z(Sp(W))$ acts via the sign character. Moreover, using the
description for the action of $-I$ given in Remark \ref{Schrodinger Model} 
\textit{(C), }it is easy to verify that for $n$ odd 
\begin{equation*}
\dim (\omega _{\psi ,1})=\left\{ 
\begin{array}{c}
\frac{q^{n}+1}{2},\text{ \ if \ }q\equiv 1\text{ }\func{mod}\text{ }4; \\ 
\frac{q^{n}-1}{2},\text{ \ if \ }q\equiv 3\text{ }\func{mod}\text{ }4,%
\end{array}%
\right. \text{ \ \ and \ \ }\dim (\omega _{\psi ,sgn})=\left\{ 
\begin{array}{c}
\frac{q^{n}-1}{2},\text{ \ if \ }q\equiv 1\text{ }\func{mod}\text{ }4; \\ 
\frac{q^{n}+1}{2},\text{ \ if \ }q\equiv 3\text{ }\func{mod}\text{ }4,%
\end{array}%
\right. \text{ }
\end{equation*}%
and for $n$ even 
\begin{equation*}
\dim (\omega _{\psi ,1})=\frac{q^{n}+1}{2}\text{ , \ \ and\ \ \ }\dim
(\omega _{\psi ,sgn})=\frac{q^{n}-1}{2}\text{. }
\end{equation*}%
Now, substituting $n=1$ we have $Sp(W)=SL_{2}(\mathbb{F}_{q})$ and using the
two possible oscillator representations (see Proposition \ref{Oscillators})
we obtain all the irreps of tensor rank one. These have dimensions $\frac{%
q\pm 1}{2}$. \ They are exactly the representations that appear "anomalous"
in the philosophy of cusp forms---see Section \ref{DD-DP}.

Next, we want information on the tensor rank two irreps.

\section{\textbf{The eta Correspondence}}

How to get information on the tensor rank two irreps of $G=SL_{2}(\mathbb{F}%
_{q})$?\textbf{\ }This section will include an answer to this question. In
fact, we introduce a systematic method called $\eta $\textit{-correspondence}
that works uniformly for all tensor ranks $k=0,1,2,$ generalizes the
construction of tensor rank one irreps described in Section \ref{HOR}, and
exists for all classical groups over finite and local fields \cite%
{Gurevich-Howe15, Gurevich-Howe17, Howe17-1, Howe17-2}.

\subsection{The $\left( O_{k\pm },SL_{2}(\mathbb{F}_{q})\right) $\textbf{\
Dual Pair\label{DP}}}

Consider the vector space $U=\mathbb{F}_{q}^{k}$ and let $\beta $ be an
inner product (i.e., a non-degenerate symmetric bilinear form) on $U$. The
pair $(U,\beta )$ is called a \textit{quadratic space.} In this note we will
be interested only in the cases $k=0,1,2.$ In these cases we have the
following examples:

\begin{itemize}
\item $k=0$: The zero space $U=0$ with $\beta =0.$

\item $k=1$: The line $U=\mathbb{F}_{q}$ with the form $\beta _{1+}(x)=x^{2}$
or with the form $\beta _{1-}(x)=\varepsilon x^{2},$ where $\varepsilon $ is
a non-square in $\mathbb{F}_{q}$.

\item $k=2$: The plane $U=\mathbb{F}_{q}^{2}$ with

\begin{itemize}
\item the form $\beta _{2+}(x,y)=x^{2}-y^{2}$. In this case---called \textit{%
hyperbolic plane---i.e., there are two lines on which the form is
identically zero. }

\item the form $\beta _{2-}(x,y)=x^{2}-\varepsilon y^{2},$ where $%
\varepsilon $ is a non-square in $\mathbb{F}_{q}$. In this case---called 
\textit{anisotropic plane}---there is no line on which the form is
identically zero.
\end{itemize}
\end{itemize}

It is well known \cite{Lam05} that up to isometry there are only two
quadratic spaces of dimension $k\neq 0$ over $\mathbb{F}_{q}.$ In
particular, for $k=0,1,2,$ the above list of examples is exhaustive.

Finally, for a quadratic space $(U,\beta )$ we denote by $O_{\beta }$ the
associated isometry group and call it the orthogonal group. We might also
denote this group by $O_{k+}$ or $O_{k-}$ according as the form $\beta $ is,
respectively, $\beta _{k+}$ or $\beta _{k-},$ or by $O_{k\pm }$ to mean
either one of these groups.

Now, let us take the plane $V=\mathbb{F}_{q}^{2}$ with the symplectic form%
\begin{equation*}
\left\langle v,v^{\prime }\right\rangle _{V}=v^{t}%
\begin{pmatrix}
0 & 1 \\ 
-1 & 0%
\end{pmatrix}%
v^{\prime },
\end{equation*}%
and consider the vector space $U\otimes V$---the tensor product of $U$ and $%
V $ \cite{Atiyah-Macdonald69}. It has a natural structure of a symplectic
space, with the symplectic form given by $\beta \otimes \left\langle
,\right\rangle _{V}.$ The groups $O_{\beta }$ and $G=SL_{2}(\mathbb{F}_{q})$
act on $U\otimes V$ via their actions on the first and second factors,
respectively,%
\begin{equation*}
O_{\beta }\curvearrowright U\otimes V\curvearrowleft G.
\end{equation*}%
Both actions preserve the form $\beta \otimes \left\langle ,\right\rangle
_{V},$ and moreover the action of $O_{\beta }$ commutes with that of $G$,
and vice versa. In particular, we have a map 
\begin{equation}
O_{\beta }\times G\longrightarrow Sp(U\otimes V),  \label{Embedding}
\end{equation}%
that embeds each of the two factors $O_{\beta }$ and $G$ in $Sp(U\otimes V),$
and the two images form a pair of commuting subgroups. In fact, each is the
full centralizer of the other inside $Sp(U\otimes V).$ Thus, the pair $%
(O_{\beta },G)$ forms what has been called in \cite{Howe73-1} a \textit{dual
pair} of subgroups of $Sp(U\otimes V)$.

\subsection{\textbf{The eta Correspondence\label{eC}}}

Let us fix for the rest of this paper a non-trivial character $\psi $ of $%
\mathbb{F}_{q}.$ Consider the symplectic space $W=U\otimes V$ with the form $%
\left\langle ,\right\rangle =\beta \otimes \left\langle ,\right\rangle _{V}$
described in Section \ref{DP} above. In this setting, we have the oscillator
representation $\omega _{U\otimes V}=\omega _{U\otimes V,\psi }$ of $%
Sp(U\otimes V)$ given by (\ref{ome-psi}).

Note that for $k=0,$ $\omega _{U\otimes V}$ is the trivial representation,
and for $k=1$ it is, respectively, the oscillator representation $\omega
_{\psi +}$ or $\omega _{\psi -}$ of $G=SL_{2}(\mathbb{F}_{q})$, according as
the form $\beta $ is $\beta _{1+}$ or $\beta _{1-}.$

\begin{proposition}
\label{omega-res-G}Assume that $\dim (U)=2$. As a representation of $%
G=SL_{2} $\bigskip $(\mathbb{F}_{q})$, in case the form $\beta $ is $\beta
_{2+},$ we have 
\begin{equation}
\omega _{U\otimes V|G}=\left\{ 
\begin{array}{c}
\omega _{\psi +}\otimes \omega _{\psi +}\text{ \ if \ }q\equiv 1\text{ }%
\func{mod}\text{ }4; \\ 
\omega _{\psi +}\otimes \omega _{\psi -}\text{ \ if \ }q\equiv 3\text{ }%
\func{mod}\text{ }4,%
\end{array}%
\right.  \label{om-res-G-2+}
\end{equation}%
and in case $\beta $ is $\beta _{2-},$ we have%
\begin{equation}
\omega _{U\otimes V|G}=\left\{ 
\begin{array}{c}
\omega _{\psi +}\otimes \omega _{\psi -}\text{ \ if \ }q\equiv 1\text{ }%
\func{mod}\text{ }4; \\ 
\omega _{\psi +}\otimes \omega _{\psi +}\text{ \ if \ }q\equiv 3\text{ }%
\func{mod}\text{ }4.%
\end{array}%
\right.  \label{om-res-G-2-}
\end{equation}
\end{proposition}

For a proof of Proposition \ref{omega-res-G} see Appendix \ref{P-omega-res-G}%
.\medskip

For the rest of this note, if not otherwise stated, we assume that the
dimension of $U$ is $k=0,1$ or $2.$ Note that, by Proposition \ref%
{omega-res-G}, every representation of $G$ that appears in $\omega
_{U\otimes V}$ is of tensor rank $\leq k$. Let us denote by $\omega
_{U\otimes V,k}$ the $k$\textit{-spectrum} of $G$ in $\omega _{U\otimes V}$,
i.e., the subspace of $\omega _{U\otimes V}$ consisting of the isotypic
components of tensor rank $k$ irreps of $G.$ The group algebras $\mathcal{A}%
_{O_{\beta }}=%
%TCIMACRO{\U{2102} }%
%BeginExpansion
\mathbb{C}
%EndExpansion
\lbrack \omega _{U\otimes V}(O_{\beta })],$ $\mathcal{A}_{G}=%
%TCIMACRO{\U{2102} }%
%BeginExpansion
\mathbb{C}
%EndExpansion
\lbrack \omega _{U\otimes V}(G)]$ act on $\omega _{U\otimes V,k}$. Let us
denote by $\widehat{G}_{k}$ the irreps of $G$ of tensor rank $k.$ The
following is a key result:

\begin{theorem}[$\protect\eta $\textbf{-correspondence}]
\label{etaC-Thm}We have

\begin{enumerate}
\item The algebras $\mathcal{A}_{O_{\beta }}$ and $\mathcal{A}_{G}$ are each
other's commutant in $End(\omega _{U\otimes V,k}).$

In particular, the double commutant theorem (see Appendix \ref{DCT}) implies
that, we have an injection 
\begin{equation}
\tau \mapsto \eta (\tau ),  \label{etaC}
\end{equation}%
from a subset of $\widehat{O_{\beta }}$ to $\widehat{G}_{k}$.

\item Using the two possible $\beta $s we obtain all of \ $\widehat{G}_{k}$,
i.e.,%
\begin{equation*}
\eta (\widehat{O_{k+}})\tbigcup \eta (\widehat{O_{k-}})=\widehat{G}_{k},
\end{equation*}%
where $\eta (\widehat{O_{k\pm }})$ denotes the image of the mapping (\ref%
{etaC}).
\end{enumerate}
\end{theorem}

We call (\ref{etaC}) the \underline{eta correspondence}.\smallskip

In \cite{Gurevich-Howe15, Gurevich-Howe17}, following initial ideas from 
\cite{Howe73-1}, we introduced techniques proving Theorem \ref{etaC-Thm} in
generality for any dual pair when one of the members is moving in a Witt
tower.

In this note we offer more elementary and explicit verification.

\subsection{\textbf{Explicit Description of the eta Correspondence\label%
{Exp-etaC}}}

We are interested in getting information on the restriction of the
oscillator representation $\omega _{U\otimes V}$ of $Sp(U\otimes V)$ to $%
G=Sp(V)=SL_{2}(\mathbb{F}_{q}).$ Consider the intertwining number $(\omega
_{U\otimes V},\omega _{U\otimes V})_{G}=\dim End_{G}(\omega _{U\otimes V})$.
In case of $\dim (U)=k=0,1,$ we know (see Proposition \ref{IO} for $k=1$)
that this number is $1$ and $2$, respectively.

\begin{proposition}
\label{IO-G}Suppose $\dim (U)=2.$ Then,%
\begin{equation}
(\omega _{U\otimes V},\omega _{U\otimes V})_{G}=2q+1.  \label{InO-G}
\end{equation}
\end{proposition}

For a proof of Proposition \ref{IO-G} see Appendix \ref{P-IO-G}.\medskip

Now, we consider the restriction, i.e., the pullback via the map (\ref%
{Embedding}), of $\omega _{U\otimes V}$ to the product $O_{\beta }\times G.$
We decompose this restriction into isotypic components for $O_{\beta }$%
\begin{equation}
\omega _{U\otimes V|O_{\beta }\times G}=\sum_{\tau \in \widehat{O_{\beta }}%
}\tau \otimes \Theta (\tau ),  \label{Res-DP}
\end{equation}%
where $\Theta (\tau )$ is a representation of $G.$

We want now to make the decomposition (\ref{Res-DP}) explicit. In most cases 
$\Theta (\tau )$ is equal to $\eta (\tau )$, however, there are cases when $%
\eta (\tau )$ is a proper sub-representation, and cases when $\Theta (\tau )$
does not contribute to the $k$-spectrum at all.

In case $\dim (U)=0$ we make the convention $O_{0}=\{1\}$ and so $\tau =1$
and $\Theta (1)=1$ the trivial representation.

\subsubsection{\textbf{The }$O_{1\pm }$-$G$ \textbf{decomposition\label%
{O1-G-decomp}}}

Let us denote by $\omega _{\psi \pm }$ the two oscillator representations $%
\omega _{U\otimes V}$ of $SL_{2}(\mathbb{F}_{q})=Sp(U\otimes V)$ in the
cases that $U$ is equipped with the forms $\beta _{1\pm }$, respectively.
The orthogonal groups $O_{1\pm }=\{\pm 1\}$ act on the representations $%
\omega _{\psi \pm },$ respectively. In these cases the isotypic components
in (\ref{Res-DP}) are $1\otimes \omega _{\psi \pm ,1}$ and $sgn\otimes
\omega _{\psi \pm ,sgn},$ where $1,$ $sgn,$ are the trivial and sign
representations of $O_{1\pm },$ and $\omega _{\psi \pm ,1},$ $\omega _{\psi
\pm ,sgn},$ are the four irreps that we discussed in Section \ref{RO}. In
particular, in these cases $\eta (\tau )=\Theta (\tau )$, and we verified
Theorem \ref{etaC-Thm} explicitly.

\subsubsection{\textbf{The }$O_{2\pm }$-$G$ \textbf{decomposition}}

We first recall some information on the orthogonal groups $O_{2\pm }$ and
their irreps, and then we describe the decomposition (\ref{Res-DP})
associated with these groups.\medskip

\paragraph{\textbf{The irreps of} $O_{2\pm }$}

The groups $O_{2\pm }$ fit into a short exact sequence 
\begin{equation}
1\rightarrow SO_{2\pm }\rightarrow O_{2\pm }\rightarrow \{\pm 1\}\rightarrow
1,  \label{ses}
\end{equation}%
where the third morphism is $\det $, and $SO_{2\pm }$ stand for the special
orthogonal groups. It might be convenient for us to choose a splitting of (%
\ref{ses}), i.e., to consider a reflection%
\begin{equation}
r\in O_{2\pm }\smallsetminus SO_{2\pm },\text{ \ \ }r^{2}=1.  \label{r}
\end{equation}%
Then $O_{2\pm }=SO_{2\pm }\rtimes <r>$ and $rsr=s^{-1}$ for every $s\in
SO_{2\pm }.$ Moreover, the groups $SO_{2\pm }$ are cyclic with $%
SO_{2+}\simeq \mathbb{F}_{q}^{\ast }$ and $SO_{2-}\simeq \{\zeta \in \mathbb{%
F}_{q^{2}}^{\ast };$ $Norm(\zeta )=1\}$ - the norm one elements in a
quadratic extension of $\mathbb{F}_{q}$. In particular, $\#(SO_{2+})=q-1$
and $\#(SO_{2-})=q+1.$%
%TCIMACRO{%
%\FRAME{fhFU}{3.7403in}{0.9591in}{0pt}{\Qcb{The irreps of $O_{2\pm }$.}}{\Qlb{irro2}}{irro2.bmp}{%
%\special{language "Scientific Word";type "GRAPHIC";maintain-aspect-ratio TRUE;display "USEDEF";valid_file "F";width 3.7403in;height 0.9591in;depth 0pt;original-width 8.2365in;original-height 2.0833in;cropleft "0";croptop "1";cropright "1";cropbottom "0";filename '../Gurevich-Howe's paper/IrrO2.bmp';file-properties "XNPEU";}} }%
%BeginExpansion
\begin{figure}[h]\centering
\includegraphics
%[natheight=2.0833in, natwidth=8.2365in, height=0.9591in, width=3.7403in]
{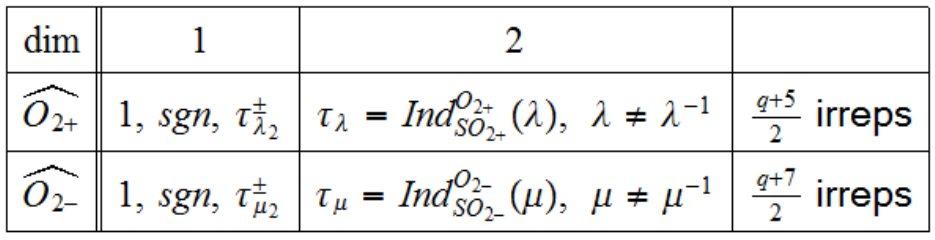}
\caption{The irreps of $O_{2\pm }$.}\label{irro2}
\end{figure}
%EndExpansion
\ 

Let us recall the irreps of $O_{2+}$---see Figure \ref{irro2} for a summary.
There are $\frac{q+5}{2}$ of them and they can be realized as follows. For
every character $\lambda $ of $SO_{2+\text{ }}$we consider the induced
representation 
\begin{equation}
\tau _{\lambda }=Ind_{SO_{2+}}^{O_{2+}}(\lambda )=\{f:O_{2+}\rightarrow 
%TCIMACRO{\U{2102} }%
%BeginExpansion
\mathbb{C}
%EndExpansion
;\text{ \ }f(sh)=\lambda (s)f(h),\text{ for }s\in SO_{2+},\text{ }h\in
O_{2+}\}.  \label{tau_lam}
\end{equation}%
Note that

\begin{itemize}
\item $\dim (Ind_{SO_{2+}}^{O_{2+}}(\lambda ))=2$ for each $\lambda .$

\item $r$ (\ref{r}) induces an isomorphism $Ind_{SO_{2+}}^{O_{2+}}(\lambda )%
\widetilde{\rightarrow }Ind_{SO_{2+}}^{O_{2+}}(\lambda ^{-1}),$ $f\mapsto
f_{r}$, where $f_{r}(h)=f(rh).$

\item $Ind_{SO_{2+}}^{O_{2+}}(\lambda )_{|SO_{2+}}=\lambda \oplus \lambda
^{-1}.$
\end{itemize}

In particular, the $\tau _{\lambda }$'s with $\lambda \neq \lambda ^{-1}$
give the $\frac{q-3}{2}$ irreps of dimension $2.$ In addition, there are $4$
irreps of dimension $1$: The trivial representation $1$, the sign character $%
sgn=\det $, and the two components $\tau _{\lambda _{2}}^{\pm }$ of $%
Ind_{SO_{2+}}^{O_{2+}}(\lambda _{2})$ on which the reflection (\ref{r}) acts%
\footnote{%
The notation $\tau _{\lambda _{2}}^{\pm }$ depends on $r.$ If $r^{\prime }$
is another reflection, then $r^{\prime }=sr$ for unique $s\in SO_{2\pm }.$}
by $+1$ or $-1$, respectively, where $\lambda _{2}$ is the unique
non-trivial character of $SO_{2+}$ such that $\lambda _{2}^{2}=1.$

The $\frac{q+7}{2}$ irreps of $O_{2-}$ are described in the same manner
starting with $\frac{q-1}{2}$ of them of dimension $2,$ induced from
characters $\mu $, $\mu \neq \mu ^{-1},$ of $SO_{2-},$ etc.---see Figure \ref%
{irro2} for the outcome in this case.\medskip

\paragraph{\textbf{The decomposition}}

\begin{figure}[h]\centering
\includegraphics
%[natheight=3.6806in, natwidth=13.2221in, height=1.7028in, width=7.3033in]
{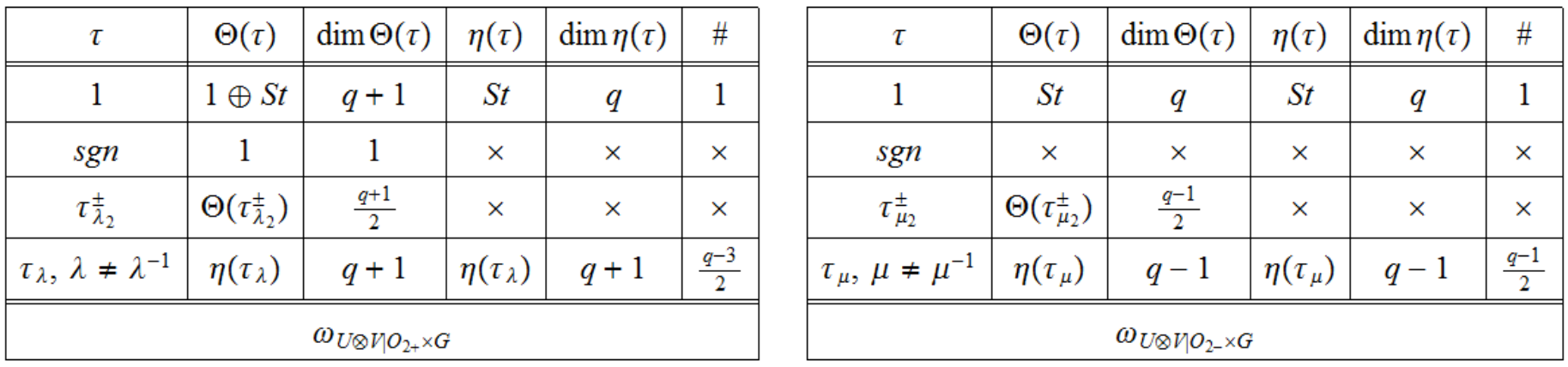}%
\caption{Explicit description of $\protect\omega _{U\otimes V|O_{2\pm
}\times G}$.}\label{omega-res-o2xg}%
\end{figure}

We need to describe (\ref{Res-DP}). Let us first treat the $O_{2-}$
case---see Figure \ref{omega-res-o2xg} (right) for a summary.%
%TCIMACRO{%
%\FRAME{fhFU}{7.3033in}{1.7028in}{0pt}{\Qcb{Explicit description of $\omega _{U\otimes V|O_{2\pm }\times G}$.}}{\Qlb{omega-res-o2xg}}{omega-res-o2xg.bmp}{%
%\special{language "Scientific Word";type "GRAPHIC";maintain-aspect-ratio TRUE;display "USEDEF";valid_file "F";width 7.3033in;height 1.7028in;depth 0pt;original-width 13.2221in;original-height 3.6806in;cropleft "0";croptop "1";cropright "1";cropbottom "0";filename '../Gurevich-Howe's paper/omega-res-O2xG.bmp';file-properties "XNPEU";}}}%
%BeginExpansion
%EndExpansion

\begin{proposition}[\textbf{The }$\protect\omega _{U\otimes V|O_{2-}\times
G} $ \textbf{decomposition}]
\label{Exp-eta-O2-}Every representation $\tau \in \widehat{O_{2-}},$ $\tau
\neq sgn$, appears in $\omega _{U\otimes V|O_{2-}},$ and in that range the
corresponding representations $\Theta (\tau )$ of $G$ (\ref{Res-DP}) are
irreducible and pairwise non-isomorphic. In particular, Part (1) of Theorem %
\ref{etaC-Thm} holds true. Moreover, for the description of $\Theta (\tau )$
and $\eta (\tau )$ we have

\begin{enumerate}
\item The representation $\Theta (1)$ is of tensor rank $2$ and dimension $%
q. $ We denote it by $\eta (1)$ or by $St$ since it is known as the
Steinberg representation.

\item The representations $\Theta (\tau _{\mu }),$ $\mu \neq \mu ^{-1}$, are
of tensor rank $2$ and dimension $q-1$. We denote them by $\eta (\tau _{\mu
})$.

\item The representations $\Theta (\tau _{\mu _{2}}^{\pm })$ are the two
tensor rank $1$ irreps of dimension $\frac{q-1}{2}.$
\end{enumerate}
\end{proposition}

For a proof of Proposition \ref{Exp-eta-O2-} see Appendix \ref%
{P-Exp-eta-O2-+}.\medskip\ 

Parts (1) and (2) of Proposition \ref{Exp-eta-O2-} give an explicit
description of the eta correspondence (\ref{etaC}).\smallskip\ 

Next, we treat the $O_{2+}$ case---see Figure \ref{omega-res-o2xg} (left)
for a summary.

\begin{proposition}[\textbf{The }$\protect\omega _{U\otimes V|O_{2+}\times
G} $ \textbf{decomposition}]
\label{Exp-eta-O2+}Every representation $\tau \in \widehat{O_{2+}}$ appears
in $\omega _{U\otimes V|O_{2+}}$. Moreover, for the description of $\Theta
(\tau )$ and $\eta (\tau ),$ we have

\begin{enumerate}
\item The representation $\Theta (1)=1\oplus St,$ where $St$ is the
Steinberg representation of tensor rank $2$ and dimension $q$. We denote it
by $\eta (1).$ In addition, $\Theta (sgn)=1$.

\item The representations $\Theta (\tau _{\lambda }),$ $\lambda \neq \lambda
^{-1},$ are irreducible, pairwise non-isomorphic, of tensor rank $2$ and
dimension $q+1$. We denote them by $\eta (\tau _{\lambda })$.

\item The representations $\Theta (\tau _{\lambda _{2}}^{\pm })$ are the two
tensor rank $1$ irreps of dimension $\frac{q+1}{2}.$

In particular, note that on the tensor rank $2$ spectrum (see Section \ref%
{eC}), i.e., on the subspace $\omega _{U\otimes V,2},$ Parts (1) and (2)
give the conclusion of Part (1) of Theorem \ref{etaC-Thm}.
\end{enumerate}
\end{proposition}

For a proof of Proposition \ref{Exp-eta-O2+} see Appendix \ref%
{P-Exp-eta-O2-+}.\medskip

Parts (1) and (2) of Proposition \ref{Exp-eta-O2+} give an explicit
description of the eta correspondence (\ref{etaC}).\smallskip\ 

\begin{remark}[\textbf{Double commutatnt property}]
Note the difference between the two decompositions appearing in Propositions %
\ref{Exp-eta-O2-} and \ref{Exp-eta-O2+}. In the first case $O_{2-}$ and $G$
generate each other's commutant when acting on the space $\omega _{U\otimes
V}$. This is not what is happening in the second case of $O_{2+}$ and $G.$
However, they do generate each other's commutant---see Figure \ref%
{omega-res-o2xg} (left)---after taking away from $\omega _{U\otimes V}$ the
isotypic component of the trivial representation of $G.$
\end{remark}

\subsubsection{\textbf{Summary\label{Summary}}}

The $\eta $-correspondence that was described explicitly in this section,
gives a way to realize $q+4$ irreps of $G$. It does it for each tensor rank $%
k=0,1,2$ separately. The realization is inside specific oscillator
representations, and in each one is described in term of irreps of
corresponding orthogonal groups. In particular, $\eta (\widehat{O_{0}})=%
\widehat{G}_{0},$ $\eta (\widehat{O_{1+}})\cup \eta (\widehat{O_{1-}})=%
\widehat{G}_{1}$, and $\eta (\widehat{O_{2+}})\cup \eta (\widehat{O_{2-}})=%
\widehat{G}_{2}$, verifying Part (2) of Theorem \ref{etaC-Thm}. Finally, by
counting $\#(\widehat{G})=q+4$ so $\widehat{G}_{0}\cup \widehat{G}_{1}\cup 
\widehat{G}_{2}=\widehat{G},$ which verifies Proposition \ref{F2R}.

\section{\textbf{Application to Uniformity of the Commutator Map}}

We go back to the problem of the uniformity of the commutator map $\left[ ,%
\right] :G\times G\rightarrow G$ described in Section \ref{E-CM} for $%
G=SL_{2}(\mathbb{F}_{q})$.

\subsection{\textbf{Statement}}

We considered the function $\mathcal{N}$ on $G$ given by (\ref{N}) $\mathcal{%
N}(g)=\#(\left[ ,\right] _{g})/\#G$, i.e., the normalized cardinality of the
fiber over $g$ of the map $\left[ ,\right] .$

We wanted to show that,

\begin{theorem}
\label{T-U-N}For $\pm I\neq g\in G$ we have 
\begin{equation}
\mathcal{N}(g)=1+O(1/q).  \label{U-N}
\end{equation}
\end{theorem}

\subsection{\textbf{Verifying Cancellations in Frobenius's Character Sum }}

Recall that Frobenius's Formula (\ref{FSum}) gives a representation
theoretic interpretation of $\mathcal{N}$ as a character sum%
\begin{equation}
\mathcal{N(}g)=1+\sum_{1\neq \rho \in \widehat{G}}\frac{\chi _{\rho }(g)}{%
\dim (\rho )}.  \label{F-S}
\end{equation}%
In Section \ref{CR-TR} we looked on (\ref{U-N}) and (\ref{F-S}) and
since---see Figure \ref{cr-sl2-q}---around $q$ of the summands in (\ref{F-S}%
) are of size $\approx 1/q$ we understood that we need to verify
cancellations in the sum on the right-hand side of (\ref{F-S}). Moreover, we
noticed that for some elements $g\in G,$ the size of the summands in (\ref%
{F-S}) seem---see Figure \ref{cr-sl2-q}---larger for irreps of tensor rank $%
1 $ then for these of tensor rank $2$. Hence, at the end of Section \ref%
{CR-TR} we proposed to verify the cancellations in (\ref{F-S}) by splitting
the sum into two partial sums, one over $\widehat{G}_{1}=$ tensor rank one
irreps and one over $\widehat{G}_{2}=$tensor rank two irreps

\begin{equation*}
\sum\limits_{1\neq \rho \in \widehat{G}}\frac{\chi _{\rho }(g)}{\dim (\rho )}%
=\sum\limits_{\rho \in \widehat{G}_{1}}\frac{\chi _{\rho }(g)}{\dim (\rho )}%
+\sum\limits_{\rho \in \widehat{G}_{2}}\frac{\chi _{\rho }(g)}{\dim (\rho )},
\end{equation*}%
and explore cancellations within each sub-sum.%
%TCIMACRO{%
%\FRAME{fhFU}{3.7273in}{1.6025in}{0pt}{\Qcb{Partial sums over irreps of tensor rank one and two for $SL_{2}(\mathbb{F}_{101}).$}}{\Qlb{s1-s2-sl2-101}}{s1-s2-sl2-101.bmp}{%
%\special{language "Scientific Word";type "GRAPHIC";maintain-aspect-ratio TRUE;display "USEDEF";valid_file "F";width 3.7273in;height 1.6025in;depth 0pt;original-width 8.1803in;original-height 2.8193in;cropleft "0";croptop "1";cropright "1";cropbottom "0";filename '../Gurevich-Howe's paper/S1-S2-SL2-101.bmp';file-properties "XNPEU";}}}%
%BeginExpansion
\begin{figure}[h]\centering
\includegraphics
%[natheight=2.8193in, natwidth=8.1803in, height=1.6025in, width=3.7273in]
{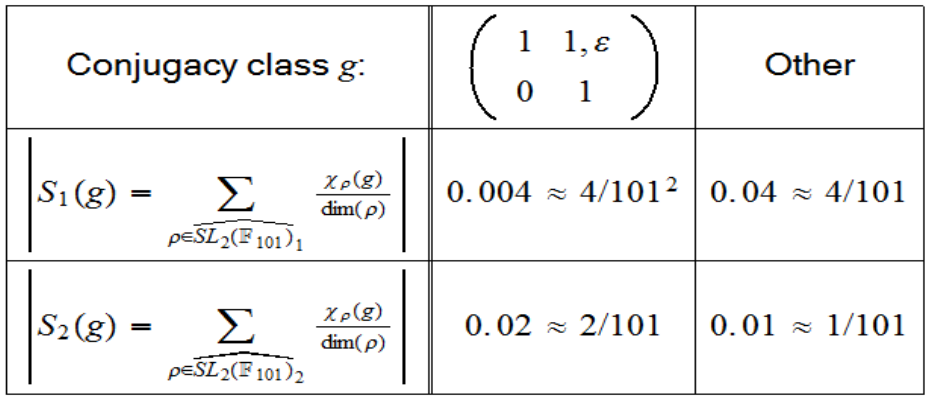}%
\caption{Partial sums over irreps of tensor rank one and two for $SL_{2}(%
\mathbb{F}_{101}).$}\label{s1-s2-sl2-101}%
\end{figure}%
%EndExpansion

The numerics appearing in Figure \ref{s1-s2-sl2-101} suggest that the above
idea might work. Indeed, using the information on $\widehat{G}_{1}$ and $%
\widehat{G}_{2}$ supplied by the $\eta $-correspondence, we can show that

\begin{proposition}
\label{S1S2}Let $\pm I\neq g\in G$. We have%
\begin{equation}
S_{1}(g)=\sum\limits_{\rho \in \widehat{G}_{1}}\frac{\chi _{\rho }(g)}{\dim
(\rho )}=\left\{ 
\begin{array}{c}
O(1/q^{2}),\text{\ \ if \ }g\sim \tbinom{1\text{ \ \ }1,\varepsilon }{0\text{
\ \ \ \ }1}\text{;} \\ 
O(1/q),\text{ \ \ \ \ \ \ \ \ \ \ other .\ \ \ }%
\end{array}%
\right.  \label{S1}
\end{equation}%
and 
\begin{equation}
S_{2}(g)=\sum\limits_{\rho \in \widehat{G}_{2}}\frac{\chi _{\rho }(g)}{\dim
(\rho )}=O(1/q).  \label{S2}
\end{equation}
\end{proposition}

Proposition \ref{S1S2} implies Theorem \ref{T-U-N}.

\subsection{\textbf{Proof of Proposition \protect\ref{S1S2}}}

The $\eta $-correspondence enables to relate the sums (\ref{S1}) and (\ref%
{S2}) that we want to estimate, to some other sums that we can actually
compute using only two inputs: The character of the oscillator
representation, for which we have explicit formulas \cite{Gerardin77,
Gurevich-Hadani07, Howe73-1, Howe73-2, Issacs73, Teruji08}; and the
character of the Steinberg representation, for which a simple description is
known \cite{Fulton-Harris91}.

Let us denote by $\chi _{\omega _{\psi \pm }}$ the characters of the two
oscillator representations of $G=SL_{2}(\mathbb{F}_{q})$ associated with a
character $1\neq \psi $ of $\mathbb{F}_{q}$ (see Section \ref{O1-G-decomp}).
The formulas (see \textit{loc. sit.}) in particular give:

\begin{fact}
For $I\neq g\in G$ we have

\begin{enumerate}
\item The oscillator characters satisfies,%
\begin{equation}
\chi _{\omega _{\psi \pm }}(g)=\left\{ 
\begin{array}{c}
\pm \gamma (\psi ),\text{ \ if \ }g\sim \tbinom{1\text{ \ \ }1}{0\text{\ \ \ 
}1}; \\ 
\text{ } \\ 
\mp \gamma (\psi ),\text{ \ if \ }g\sim \tbinom{1\text{ \ \ }\varepsilon }{0%
\text{\ \ \ }1}; \\ 
\text{\ \ } \\ 
\QTOVERD( ) {-\det (g-I)}{q},\ \ \ \ \text{other, \ \ }%
\end{array}%
\right. \text{\ \ }  \label{X_omega}
\end{equation}%
where $\gamma (\psi )=\tsum\limits_{z\in \mathbb{F}_{q}}\psi (\frac{1}{2}%
z^{2})$ is the usual quadratic Gauss sum \cite{Ireland-Rosen90}, $\QTOVERD(
) {\cdot }{q}$ is the Legendre symbol, and $\sim $ means in the same
conjugacy class.

\item The Steinberg Character satisfies,%
\begin{equation}
\chi _{St}(g)=\#(\mathcal{L}^{g})-1,  \label{X_St}
\end{equation}%
where $\mathcal{L}^{g}$ is the collection of lines in the plane $V$ that are
fixed by $g$.
\end{enumerate}
\end{fact}

\underline{\textbf{The tensor rank one sum}}\smallskip

First, we estimate the sum $S_{1}(g)$ (\ref{S1}) for $g\neq \pm I$.

Let us write $\chi _{\omega _{\psi +}}=\chi _{++}+\chi _{+-},$ where $\chi
_{++}$ and $\chi _{+-}$ denote the characters of the components of $\omega
_{\psi +}$ of dimensions $\frac{q+1}{2}$ and $\frac{q-1}{2},$ respectively,
and likewise $\chi _{\omega _{\psi -}}=\chi _{-+}+\chi _{--}.$

We have%
\begin{eqnarray}
S_{1} &=&\frac{\chi _{+-}}{\frac{q-1}{2}}+\frac{\chi _{--}}{\frac{q-1}{2}}+%
\frac{\chi _{++}}{\frac{q+1}{2}}+\frac{\chi _{-+}}{\frac{q+1}{2}}
\label{S_1} \\
&=&(\frac{1}{\frac{q-1}{2}}-\frac{1}{\frac{q+1}{2}})(\chi _{+-}+\chi _{--})+%
\frac{1}{\frac{q+1}{2}}(\chi _{++}+\chi _{-+}+\chi _{+-}+\chi _{--}).  \notag
\end{eqnarray}%
But, Formula (\ref{X_omega}) implies that 
\begin{eqnarray}
\chi _{++}(g)+\chi _{-+}(g)+\chi _{+-}(g)+\chi _{--}(g) &=&\chi _{\omega
_{\psi +}}(g)+\chi _{\omega _{\psi -}}(g)  \label{X++-} \\
&=&\left\{ 
\begin{array}{c}
0\text{, \ \ if \ }g\sim \tbinom{1\text{ \ \ }1,\varepsilon }{0\text{\ \ \ \ 
}1\text{ }};\text{ } \\ 
\pm 2\text{, \ \ \ \ \ \ \ \ \ other},\text{ \ \ \ \ \ \ \ }%
\end{array}%
\right.  \notag
\end{eqnarray}

and, moreover, from the definition of $\chi _{+-}$ and $\chi _{--}$ and
Formula (\ref{X_omega}), we get%
\begin{equation}
\chi _{+-}(g)+\chi _{--}(g)=\left\{ 
\begin{array}{c}
-1\text{, \ \ \ \ \ \ \ if \ }g\sim \tbinom{\pm 1\text{ \ \ }1,\varepsilon }{%
\text{ \ }0\text{\ \ \ \ }\pm 1\text{ }}; \\ 
\left\vert \cdot \right\vert \leq 2\text{, \ \ \ \ \ \ \ \ \ \ other}.\text{
\ \ \ \ \ \ }%
\end{array}%
\right.  \label{X_q-1}
\end{equation}%
So overall, by a combination (\ref{S_1}), (\ref{X++-}), and (\ref{X_q-1}),
we have%
\begin{equation*}
\left\vert S_{1}(g)\right\vert =\left\{ 
\begin{array}{c}
\text{ \ }\frac{4}{q^{2}-1}\text{, \ \ \ \ \ \ \ if \ }g\sim \tbinom{1\text{
\ \ }1,\varepsilon }{0\text{\ \ \ \ }1\text{ }}; \\ 
\frac{4}{q+1}+\delta \text{, \ \ \ \ \ \ \ \ \ \ other, \ \ \ \ }%
\end{array}%
\right.
\end{equation*}%
where $\left\vert \delta \right\vert \leq \frac{8}{q^{2}-1}.$

This completes a quantitative form of (\ref{S1}).\smallskip

\underline{\textbf{The tensor rank two sum}}\smallskip

Next, we estimate the sum $S_{2}(g)$ (\ref{S2}) for $g\neq \pm 1$. It has a
similar structure to $S_{1}(g).$

Using the explicit description of the $\eta $-correspondence---see Figure %
\ref{omega-res-o2xg}---we have%
\begin{eqnarray}
S_{2} &=&\tsum\limits_{\tau _{\mu }\in \widehat{O_{2-}},\text{ }\mu \neq \mu
^{-1}}\frac{\chi _{\eta (\tau _{\mu })}}{q-1}+\tsum\limits_{\tau _{\lambda
}\in \widehat{O_{2+}},\text{ }\lambda \neq \lambda ^{-1}}\frac{\chi _{\eta
(\tau _{\lambda })}}{q+1}+\frac{\chi _{St}}{q}  \label{S_2} \\
&=&(\frac{1}{q-1}-\frac{1}{q+1})\tsum\limits_{\tau _{\mu }\in \widehat{O_{2-}%
},\text{ }\mu \neq \mu ^{-1}}\chi _{\eta (\tau _{\mu })}  \notag \\
&&+\frac{1}{q+1}(\text{ }\underset{\mathcal{X}_{2-}}{\underbrace{%
\tsum\limits_{\tau _{\mu }\in \widehat{O_{2-}},\text{ }\mu \neq \mu
^{-1}}\chi _{\eta (\tau _{\mu })}}}+\underset{\mathcal{X}_{2+}}{\underbrace{%
\tsum\limits_{\tau _{\lambda }\in \widehat{O_{2+}},\text{ }\lambda \neq
\lambda ^{-1}}\chi _{\eta (\tau _{\lambda })}}}+\chi _{St}\text{ })+(\frac{1%
}{q}-\frac{1}{q+1})\chi _{St}.  \notag
\end{eqnarray}

We would like to estimate the various terms in (\ref{S_2}).

Consider the sums $\mathcal{X}_{2-}$ and $\mathcal{X}_{2+}$ above. Denote by 
$U_{2-}$ and $U_{2+},$ respectively, the anisotropic and hyperbolic planes,
and by $\chi _{\omega _{U_{2\pm }\otimes V}},$ respectively, the characters
of the oscillator representations of the groups $Sp(U_{2\pm }\otimes V)$
(see Section \ref{eC}). Then,%
\begin{eqnarray*}
&&2\mathcal{X}_{2-}(g)+\chi _{+-}(g)+\chi _{--}(g)+\chi _{St}(g)+2\mathcal{X}%
_{2+}(g)+\chi _{++}(g)+\chi _{-+}(g)+\chi _{St}(g)+2 \\
&=&\chi _{\omega _{U_{2+}\otimes V}}(I\otimes g)+\chi _{\omega
_{U_{2-}\otimes V}}(I\otimes g)=\left\{ 
\begin{array}{c}
0\text{, \ \ if \ }g\sim \tbinom{1\text{ \ \ }1,\varepsilon }{0\text{\ \ \ \ 
}1\text{ }};\text{ } \\ 
\pm 2\text{, \ \ \ \ \ \ \ \ \ other},\text{ \ \ \ \ \ \ \ }%
\end{array}%
\right.
\end{eqnarray*}%
\ where the first equality follows from the explicit descriptions of $\omega
_{U\otimes V|O_{2\pm }\times G}$ appearing in Figure \ref{omega-res-o2xg},
and the second equality follows from a combination of Proposition \ref%
{omega-res-G} and Formula (\ref{X_omega}). So, together with (\ref{X++-}) we
get, 
\begin{equation}
\mathcal{X}_{2-}(g)+\mathcal{X}_{2+}(g)+\chi _{St}(g)=\left\{ 
\begin{array}{c}
-1\text{, \ \ if \ }g\sim \tbinom{1\text{ \ \ }1,\varepsilon }{0\text{\ \ \
\ }1\text{ }};\text{ } \\ 
\pm 2\text{, \ \ \ \ \ \ \ \ \ other},\text{ \ \ \ \ \ }%
\end{array}%
\right.  \label{X_2-_2+_St}
\end{equation}

Next, concerning the term $\mathcal{X}_{2-}$ in (\ref{S_2}). From the
explicit description of $\omega _{U\otimes V|O_{2-}\times G}$ appearing in
Figure \ref{omega-res-o2xg}, we obtain%
\begin{equation}
\mathcal{X}_{2-}(g)=\frac{\chi _{\omega _{U_{2-}\otimes V}}(I\otimes g)-\chi
_{+-}(g)-\chi _{--}(g)-\chi _{St}(g)}{2}.  \label{X2_}
\end{equation}%
So, using Formulas (\ref{om-res-G-2-}), (\ref{X_omega}), and (\ref{X_q-1}), (%
\ref{X_St}), inserted in (\ref{X2_}), we get 
\begin{equation}
\mathcal{X}_{2-}(g)=\left\{ 
\begin{array}{c}
-\frac{q-1}{2}\text{, \ \ \ \ \ \ \ if \ }g\sim \tbinom{1\text{ \ \ }%
1,\varepsilon }{0\text{\ \ \ \ }1\text{ }}; \\ 
\left\vert \cdot \right\vert \leq 2\text{, \ \ \ \ \ \ \ \ \ \ \ other. \ \
\ \ \ \ \ }%
\end{array}%
\right.  \label{X2-}
\end{equation}

Combining (\ref{S_2}), (\ref{X_St}), (\ref{X_2-_2+_St}), and (\ref{X2-}), we
have for every $\pm I\neq g\in G,$

\begin{equation*}
\left\vert S_{2}(g)\right\vert \leq \frac{2}{q+1}+\delta ,
\end{equation*}%
\newline
where $\left\vert \delta \right\vert \leq \frac{5}{q^{2}-1}$.

This is a quantitative form of (\ref{S2}).

This completes effective proofs (with explicit bounds) for Proposition \ref%
{S1S2} and Theorem \ref{T-U-N}.

\appendix

\section{\textbf{The Weyl transform}}

Let us denote by $\mathcal{H}_{\psi }$ a vector space supporting the
Heisenberg representation $\pi _{\psi }$ of the Heisenberg group $H=H(W)$
defined in Section \ref{RHG}. Various calculations done in this note involve
the space $End(\mathcal{H}_{\psi })$ of \textit{all} linear transformation
on $\mathcal{H}_{\psi }.$ Interestingly, this operator space has an
invariant description \cite{Gurevich-Hadani07, Howe73-1} that was discovered
by Weyl in the context of quantum mechanics \cite{Weyl27}.

\subsection{\textbf{Definition of the Weyl Transform}}

Consider the vector space $L^{2}(W)$ of complex valued functions on the
symplectic space $W.$ Then, we have a canonical isomorphism, called the 
\textit{Weyl transform,} 
\begin{equation}
\left\{ 
\begin{array}{c}
\mathcal{W}:End(\mathcal{H}_{\psi })\widetilde{\longrightarrow }L^{2}(W) \\ 
\mathcal{W}(T)(w)=\frac{1}{\dim (\mathcal{H}_{\psi })}trace(T\circ \pi
_{\psi }(-w)),%
\end{array}%
\right.  \label{WT}
\end{equation}%
for every $T\in End(\mathcal{H}_{\psi }),$ $w\in W.$

Here is a useful application of the Weyl transform (\ref{WT}).

\subsection{\textbf{Application for Intertwining Numbers}}

We describe application for the computation of intertwining numbers for
restrictions of the oscillator representation $\omega _{\psi }$ to subgroups 
$K<Sp(W).$ Note that the identification $\mathcal{W}$ intertwines the
conjugation action of $Sp(W)$ on $End(\mathcal{H}_{\psi })$ with its
permutation action on $L^{2}(W).$ Consider the intertwining number $\left(
\omega _{\psi },\omega _{\psi }\right) _{K}=\dim End_{K}(\mathcal{H}_{\psi
}),$ i.e., the dimension of the space of operators on $\mathcal{H}_{\psi }$
that commute with the action of $K.$ Denote by $W/K$ the set of orbits for
the action of $K$ on $W.$ It follows

\begin{corollary}
\label{IO-K}We have $\left( \omega _{\psi },\omega _{\psi }\right)
_{K}=\#(W/K).$
\end{corollary}

\section{\textbf{The Double Commutant Theorem\label{DCT}}}

In Section \ref{eC} we used the \textit{double commutant theorem }\cite%
{Weyl46} to define the $\eta $-correspondence.

\subsection{\textbf{Formulation}}

For the convenience of the reader, here is the statement.

\begin{theorem}[\textbf{Double commutant theorem}]
Let $\mathcal{H}$ be a finite dimensional vector space. Let $\mathcal{A},%
\mathcal{\mathcal{A}}^{\prime }$ $\subset End(\mathcal{H})$ be two
sub-algebras, such that

\begin{enumerate}
\item The algebra $\mathcal{A}$ acts semi-simply on $\mathcal{H}.$

\item Each of $\mathcal{A}$ and $\mathcal{\mathcal{A}}^{\prime }$ is the
full commutant of the other in $End(\mathcal{H}).$
\end{enumerate}

Then $\mathcal{A}^{\prime }$ acts semi-simply on $\mathcal{H},$ and as a
representation of\ $\mathcal{A\otimes \mathcal{A}}^{\prime }$ we have%
\begin{equation}
\mathcal{H}=\tbigoplus\limits_{i\in I}\mathcal{H}_{i}\otimes \mathcal{H}%
_{i}^{\prime },  \label{D}
\end{equation}%
where $\mathcal{H}_{i}$ are all the irreducible representations of $\mathcal{%
A}$, and $\mathcal{H}_{i}^{\prime }$ are all the irreducible representations
of $\mathcal{A}^{\prime }.$ In particular, we have a bijection between
irreducible representations of $\mathcal{A}$ and $\mathcal{A}^{\prime },$
and moreover, every isotypic component for $\mathcal{A}$ is an irreducible
representation of $\mathcal{A\otimes \mathcal{A}}^{\prime }.$

On the other hand, if $\mathcal{A},\mathcal{\mathcal{A}}^{\prime }$ commute
and as a representation of\ $\mathcal{A\otimes \mathcal{A}}^{\prime }$ \ the
decomposition (\ref{D}) holds, then each of $\mathcal{A}$ and $\mathcal{%
\mathcal{A}}^{\prime }$ is the full commutant of the other in $End(\mathcal{H%
}).$
\end{theorem}

\section{\textbf{Proofs}}

\subsection{\textbf{Proof of Proposition }\protect\ref{IO}\label{P-IO}}

\begin{proof}
Take $K=Sp(W)$ in Corollary \ref{IO-K} and note $\#(W/Sp(W))=2,$ as claimed.
\end{proof}

\subsection{\textbf{Proof of Proposition }\protect\ref{omega-res-G}\label%
{P-omega-res-G}}

\begin{proof}
Let us denote by $\mathcal{\omega }(W,\left\langle ,\right\rangle ,\psi )$
the oscillator representation associated with a symplectic vector space $%
(W,\left\langle ,\right\rangle )$ and a central character $\psi .$ It is
well known that $\mathcal{\omega }(W,a\left\langle ,\right\rangle ,\psi
)\simeq $ $\mathcal{\omega }(W,\left\langle ,\right\rangle ,\psi _{a}),$ and 
$\mathcal{\omega }((W_{1},\left\langle ,\right\rangle _{1})\oplus
(W_{2},\left\langle ,\right\rangle _{2}),\psi )\simeq \mathcal{\omega }%
(W_{1},\left\langle ,\right\rangle _{1},\psi )\otimes \mathcal{\omega }%
(W_{2},\left\langle ,\right\rangle _{2},\psi ),$ for every $a\in \mathbb{F}%
_{q}^{\ast },$ and two symplectic spaces $(W_{i},\left\langle ,\right\rangle
_{i})$, $i=1,2.$

Next, in the case $U$ endowed with the form $\beta _{2+}$ we can identify $%
(U\otimes V,\beta _{2+}\otimes \left\langle ,\right\rangle _{V}),$ as a $%
G=Sp(V)$-space, with $(V,\left\langle ,\right\rangle _{V})\oplus
(V,-\left\langle ,\right\rangle _{V}),$ and in the case $U$ is with form $%
\beta _{2-}$ we can identify $(U\otimes V,\beta _{2-}\otimes \left\langle
,\right\rangle _{V}),$ as a $G=Sp(V)$-space, with $(V,\left\langle
,\right\rangle _{V})\oplus (V,-\varepsilon \left\langle ,\right\rangle _{V})$
where $\varepsilon \in \mathbb{F}_{q}^{\ast }$ a non-square.

Now, combining the above functorial properties with Proposition \ref%
{Oscillators}, we get identities (\ref{om-res-G-2+}) and (\ref{om-res-G-2-}%
), as claimed.
\end{proof}

\subsection{Proof of Proposition \protect\ref{IO-G}\label{P-IO-G}}

\begin{proof}
We take $W=U\otimes V$ and $G=Sp(V)<Sp(W)$ in Corollary \ref{IO-K}, and
compute $\#(W/G).$ We identify the action of $G$ on $W$ with its diagonal
action on $V\times V$, and find that the orbits are

\begin{itemize}
\item For each $0\neq a\in \mathbb{F}_{q}$ the orbit $\{(u,v)\in V\times V;$ 
$\left\langle u,v\right\rangle _{V}=a\}.$

\item For each $0\neq b\in \mathbb{F}_{q}$ the orbit $\{(u,v)\in
(V\smallsetminus 0)\times (V\smallsetminus 0);$ $u=bv\}.$

\item The orbits $\{(0,v);$ $0\neq v\in V\},$ $\{(v,0);$ $0\neq v\in V\},$ $%
\{(0,0)\}.$
\end{itemize}

So, overall we have $2q+1$ orbits, as claimed.
\end{proof}

\subsection{\textbf{Proof of Propositions \protect\ref{Exp-eta-O2-} and 
\protect\ref{Exp-eta-O2+}}\label{P-Exp-eta-O2-+}}

\begin{proof}
We give a proof that works more or less in a uniform manner for both groups $%
O_{2\pm }$.

The verification will use a realization of $\omega _{U\otimes V}$ which is
convenient for the groups $O_{2\pm }.$ Indeed, $\omega _{U\otimes V}$ can be
realized on the space $L^{2}(U)$ in such a way that for $O_{2\pm }$ it is
its permutation representation. To arrive to this realization, write $V=%
\mathbb{F}_{q}\oplus \mathbb{F}_{q}$ and consider the associated Lagrangian
decomposition $U\otimes V=(U\otimes \mathbb{F}_{q})\oplus (U\otimes \mathbb{F%
}_{q})=U\oplus U$. Then we have the Schr\"{o}dinger model (see Remark \ref%
{Schrodinger Model}) realizing $\omega _{U\otimes V}$ as described above.

Recall that by Witt's theorem \cite{Lam05} the groups $O_{2\pm }$ act
transitively on the each of the sets $U_{a}=\{u\in U;$ \ $\beta _{2\pm
}(u)=a\}$ of vectors with a given value $a\in \mathbb{F}_{q}$ for the inner
product. To see what representations of $O_{2\pm }$ appear in $L^{2}(U)$ we
decompose%
\begin{equation*}
L^{2}(U)=\tbigoplus\limits_{a\in \mathbb{F}_{q}}L^{2}(U_{a}),
\end{equation*}%
and study the spaces $L^{2}(U_{a})$.

Let us start with the $q-1$ non-isotropic orbits, i.e., $U_{a}$, $a\neq 0$.
The stabilizer subgroup of a non-isotropic vector $u\in U_{a}$ is the
orthogonal group of the orthogonal complement of the line spanned by $u$,
i.e., the group generated by the reflection $r_{u}$ with respect to that
line. In particular, the groups $SO_{2\pm }$ act simply transitively on $%
U_{a}$, $a\neq 0.$

It follows that the $\det =sgn$ character can not appear in $L^{2}(U_{a}),$ $%
a\neq 0,$ since it would have to be represented by the constant function $1$
by transitivity of $SO_{2\pm }$, but this transforms by the identity
character of $O_{2\pm }$.

For the characters of $O_{2\pm }$ that are non-trivial on $SO_{2\pm }$, one
will take value $1$ on a given reflection, and one will take value $-1$.
Clearly, only the one that takes value $1$ on the reflection $r_{u}$ that
stabilizes the vector $u\in U_{a}$ can live on the orbit. That means that
each of these characters live on half of the orbits of non-isotropic vectors.

Finally, it is easy to see that each of the two-dimensional irreps of $%
O_{2\pm }$ appear once in $L^{2}(U_{a})$, $a\neq 0$.

For the double commutant property for $O_{2-},$ we know from the above
analysis that every irrep but $\det =sgn$ appears in $L^{2}(U)$, so by
Plancherel \cite{Serre77} it gives $\#(O_{2-})-1=2(q+1)-1=2q+1$ linearly
independent operators in the the commutant $End_{G}(\omega _{U\otimes V}),$
so by (\ref{InO-G}) it gives the all commutant.

In summary, the isotypic components in $\omega _{U\otimes V}$ of irreps of $%
O_{2-}$ are

\begin{itemize}
\item $1\otimes \Theta (1)$, $\Theta (1)\in \widehat{G}$ with $\dim \Theta
(1)=q$ (note that $U_{0}=\{0\}$ for $\beta _{2-}$), so $\Theta (1)=St.$

\item $\tau _{\mu }\otimes \Theta (\tau _{\mu }),$ $\mu \neq \mu ^{-1},$ $%
\Theta (\tau _{\mu })\in \widehat{G}$ pairwise non-isomorphic with $\dim
\Theta (\tau _{\mu })=q-1.$

\item $\tau _{\mu _{2}}^{\pm }\otimes \Theta (\tau _{\mu _{2}}^{\pm }),$ $%
\Theta (\tau _{\mu _{2}}^{\pm })\in \widehat{G},$ two non-isomorphic irreps
with $\dim \Theta (\tau _{\mu _{2}}^{\pm })=\frac{q-1}{2}.$
\end{itemize}

This completes the Proof of Proposition \ref{Exp-eta-O2-}.

Consider the group $O_{2+},$ note that it acts simply transitively on the
non-zero isotropic vectors, so the representation of $O_{2+}$ on $U_{0}$ is
the direct sum of the trivial representation and the regular representation.
So $L^{2}(U_{0})$ contains the $\det =sgn$ representation once, the trivial
representation twice, the two characters of $O_{2+}$ which are non-trivial
on $SO_{2+}$ each once, and each of the two-dimensional irreps of $O_{2+}$
twice.

For the double commutant property for $O_{2+}$, we have all irreps
appearing, but they only give $\#(O_{2+})=2(q-1)$ operators, so by (\ref%
{InO-G}) we are missing $3.$

We will show that

\begin{lemma}
\label{1-G}The isotypic component in $\omega _{U\otimes V}$ of the trivial
representation of $G$ is two-dimensional on which $O_{2+}$ acts by $%
sgn\oplus 1$.
\end{lemma}

This means that $O_{2+}$ only supplies two of the four operators in the
commuting algebra for this $G$ component. So if we leave out the trivial
representation of $G$, then we need only $2q+1-4=2q-3$ operators. On the
other hand, from Lemma \ref{1-G} and the analysis we have already done above
we see that, on the complement of the trivial representation of $G$ we get
all the representations of $O_{2+}$ except $sgn$, so this supplies $%
2(q-1)-1=2q-3$ operators, which is the needed number.

In summary, this will show that the isotypic components in $\omega
_{U\otimes V}$ of irreps of $O_{2+}$ are

\begin{itemize}
\item $1\otimes (1\oplus St)$.

\item $sgn\otimes 1.$

\item $\tau _{\lambda }\otimes \Theta (\tau _{\lambda }),$ $\lambda \neq
\lambda ^{-1},$ $\Theta (\tau _{\lambda })\in \widehat{G}$ pairwise
non-isomorphic with $\dim \Theta (\tau _{\lambda })=q+1.$

\item $\tau _{\lambda _{2}}^{\pm }\otimes \Theta (\tau _{\lambda _{2}}^{\pm
}),$ $\Theta (\tau _{\lambda _{2}}^{\pm })\in \widehat{G},$ two
non-isomorphic irreps with $\dim \Theta (\tau _{\lambda _{2}}^{\pm })=\frac{%
q+1}{2}.$
\end{itemize}

To prove Lemma \ref{1-G} it is enough, due to our knowledge on the size of
the commutant, to show that $sgn\oplus 1$ sits inside the $G$-invariant
subspace of $\omega _{U\otimes V}.$ To show this, we will use another
realization of $\omega _{U\otimes V}$, that helps to visualize the action of 
$G.$ Decompose $U$ into a direct sum $U=\ell \oplus \ell ^{\ast }$ of
isotropic lines. This induces the Lagrangian decomposition $U\otimes V=(\ell
\otimes V)\oplus (\ell ^{\ast }\otimes V)$, that we can further identify
with $V\oplus V$. In particular, we can realize $\omega _{U\otimes V}$ on
the space $L^{2}(V),$ where

\begin{enumerate}
\item The action of $G=Sp(V)$ is its permutation representation.

\item The action of $SO_{2+}$ is identified with the scaling action of $%
\mathbb{F}_{q}^{\ast }.$

\item The action of the reflection $r:\ell \rightarrow \ell ^{\ast }$ is via
the symplectic Fourier transform $F$ given by%
\begin{equation*}
F[f](v)=\frac{1}{q}\sum_{v^{\prime }\in V}\psi (\left\langle v,v^{\prime
}\right\rangle )f(v^{\prime }),
\end{equation*}%
for every $f\in L^{2}(V)$, $v\in V.$
\end{enumerate}

Now, denote by $\delta _{0}$ the Dirac delta function at the origin of $V$,
and by $1_{V}$ the constant function $1_{V}(v)=1$ for every $v\in V.$
Consider the functions $f_{\pm }=\delta _{0}\pm \frac{1}{q}1_{V}\in $ $%
L^{2}(V).$ Then, $f_{\pm }$ are $G$-invariant and $F(f_{\pm })=\pm f_{\pm }.$
This completes the proof of Lemma \ref{1-G} and Proposition \ref{Exp-eta-O2+}%
.
\end{proof}

\end{document}